\documentclass[11pt,a4paper, oneside]{article}

\usepackage{amsmath}
\usepackage{amsfonts}
\usepackage{amssymb}
\usepackage{amscd}
\usepackage{xcolor}
\usepackage{latexsym}
\usepackage{graphicx}
\usepackage{mathrsfs}
\usepackage{amsthm}
\usepackage{booktabs}
\usepackage{tabularx}
\usepackage{tikz}
\usetikzlibrary{arrows.meta, positioning, backgrounds, fit, shadows}
\usepackage{comment}
\usepackage{hyperref}
\usepackage{graphics}

\def\A{{\mathcal A}}

\def\C{{\mathbb C}}
\def\e{\mathrm{e}}

\def\H{{\mathcal H}}

\def\N{{\mathbb N}}
\def\NN{{\mathcal N}}

\def\Z{{\mathbb Z}}

\def\KRn{K_{\!R_n}}
\def\GRn{\Gamma_{\!\,\!R_n}}
\def\Cmax{C_{\rm max}}

\def\0{\mathfrak r}

\def\2{\mathfrak{int}}

\newtheorem{theorem}{Theorem}[section]
\newtheorem{prop}[theorem]{Proposition}

\newtheorem{remark}[theorem]{Remark}
\newtheorem{example}[theorem]{Example}
\theoremstyle{definition}

\allowdisplaybreaks

\begin{document}

\title{Elementary cellular automata: spectral properties of Koopman operator for finite periodic realizations}

\author{Q. Zhu${}^1$, S. Richard${}^2$\footnote{Supported by JSPS Grant-in-Aid for scientific research C no 26K06826}}

\date{\today}

\maketitle 

\vspace{-1cm}

\begin{quote}
\begin{itemize}
\item[1] Graduate School of Mathematics, Nagoya University, Furo-cho, Chikusa-ku, Nagoya 464-8602, Japan
\item[2] Institute for Liberal Arts and Sciences \& Graduate School of Mathematics, Nagoya University, Furo-cho, Chikusa-ku, Nagoya 464-8601, Japan
\item[] E-mail:  zhu.qingying.b0@s.mail.nagoya-u.ac.jp \\ 
${}\qquad \quad\ {}$ richard@math.nagoya-u.ac.jp
\end{itemize}
\end{quote}

\begin{abstract}
In this paper we study Koopman operators
associated with finite periodic realizations of elementary cellular automata. 
The spectra of these operators encompass various properties
of the underlying dynamical systems, independently  of any initial condition.
A characterization of these spectra is proposed, leading to a spectral
classification of all elementary cellular automata. 
This work is thought as a preliminary step for studying the spectra of Koopman operators
for infinite realizations of elementary cellular automata.
\end{abstract}

\noindent{\bf Keywords:} Elementary cellular automata, Koopman operator, finite periodic realizations, spectral classification


\section{Introduction}\label{sec_intro}

Since the seminal works of Wolfram in the $80$'s 
\cite{wolfram1983statistical,Wolfram1984,W}, 
Elementary Cellular Automata (ECA) have been extensively studied
in different settings and through multiple approaches.
These investigations led to numerous results and to various classifications
which have recently been summarized for example in \cite{BW,Martinez, SS2012, Stricker, Vispoel},
see also the references therein.
Clearly, adding something new to this extended corpus of knowledge looks like an impossible task, 
nevertheless it is the humble aim of this manuscript.

Koopman operator has been introduced a long time ago \cite{Koopman}, 
but it is only over the last couple of decades that its usefulness and importance have
been revealed again. An account on these recent developments can be found 
for example in the book \cite{MMS} or the review paper \cite{BBKK}. 
Koopman approach provides a linear framework for 
studying nonlinear dynamical systems, with the price of often ending in an infinite 
dimensional space of observables. Our investigations on ECA will precisely be based
on Koopman theory, as explained in the
rest of this introduction. Note that our interest in this setting has been triggered
by the two papers \cite{Taga1, Taga2} in which ECA are also analyzed through a Koopman
operator's approach.

What we provide in this paper is relatively simple: we study the spectral theory of  
Koopman operators for finite periodic realizations of all ECA (or more precisely of the
88 independent clusters of ECA, as introduced in \cite{li1990structure}).
For a periodic realization on $n$ cells, Koopman operator corresponds simply 
to a $2^n\times 2^n$
matrix which can be described in terms of a directed graph made of $2^n$ vertices
and $2^n$ edges. 
As explained in Section \ref{sec_setting}
this matrix contains all information on the dynamical system defined
by the elementary cellular automata, and in particular the evolution of arbitrary initial conditions.
Consequently, the spectrum of this matrix, namely the set of its eigenvalues, is very much
related to the evolution of the dynamical system. This spectrum clearly
does not contain as much information as the original matrix, since eigenvectors or generalized eigenvectors are disregarded. 
On the other hand, it can be easily evaluated in terms of the structure
of the already mentioned directed graph. In particular, there is no need to solve a 
$2^n$-dimensional characteristic polynomial or to fully analyze a huge matrix for getting
this set of eigenvalues.

Once the spectra of Koopman operators of the 88 clusters of ECA have been determined
for $n$ in a sufficiently large set, the taxonomy work can start. Indeed, instead of classifying
all ECA based on the evolution of specific initial conditions, as performed in 
the seminal work of Wolfram, they can be classified according to specific properties
of the spectra of Koopman operators.
This approach is what we call a \emph{spectral classification} of all ECA.
By definition, it does not depend on any initial condition, and does not
rely on any fixed $n$. 
In the classification proposed in Section \ref{sec_classification} 
we also tried to remove all subjective factors, 
but are aware that the initial choice of the periodic boundary condition plays a role.
This classification has also been established by keeping its future extension to infinite realizations
in mind. We do not address this extension in the present manuscript. 

Given a fixed ECA, the study of the spectrum of its Koopman operator for various $n$ also reveals 
some unexpected dependence on this parameter. From different perspectives, 
most of these properties have already been observed and are documented in the literature.
However we also mention in Section \ref{sec_additional} a few additional 
observations that we 
could not track down in the literature. For now, these are only observations and 
additional work would be necessary to turn them into proper mathematical statements.
We leave this for future work.

In summary, Koopman's approach provides a unifying framework for studying dynamical
systems in a linear setting. This allows us to investigate finite periodic realizations of elementary cellular 
automata though the eigenvalues of the corresponding Koopman operators. 
By systematically analyzing these spectra, a classification and additional structural information
are presented and discussed. These results pave the way for further  investigations on specific ECA
or for studying their infinite realization on $\Z$.

\section{Setting}\label{sec_setting}

In this section, we introduce the setting and the notations. Namely, we
firstly provide a clear definition of elementary cellular automata and of their finite 
periodic realizations.
We then associate a directed graph to each elementary cellular automaton, and define
the state transition network. Finally, we introduce Koopman operator, and
explain how its spectrum can be deduced from the state transition network.

\subsection{Elementary cellular automata and finite realizations}

\begin{figure}[htbp]
\centering
\begin{tikzpicture}[
    >=Stealth,
    bit0/.style={
        draw=black!70,
        fill=white,
        minimum width=0.48cm,
        minimum height=0.48cm,
        inner sep=0pt,
        font=\normalsize
    },
    bit1/.style={
        draw=black!80,
        fill=black,
        text=white,
        minimum width=0.48cm,
        minimum height=0.48cm,
        inner sep=0pt,
        font=\normalsize
    },
    output0/.style={
        draw=black!70,
        fill=white,
        minimum width=0.48cm,
        minimum height=0.48cm,
        inner sep=0pt,
        font=\normalsize
    },
    output1/.style={
        draw=black!80,
        fill=black,
        text=white,
        minimum width=0.48cm,
        minimum height=0.48cm,
        inner sep=0pt,
        font=\normalsize
    }
]
\node[bit1] (a111) at (0.00,0) {1};
\node[bit1] (b111) at (0.48,0) {1};
\node[bit1] (c111) at (0.96,0) {1};
\node[output0] (o111) at (0.48,-1.15) {0};

\draw[->,thick]
(0.48,-0.28) -- (o111.north);

\node[bit1] (a110) at (1.90,0) {1};
\node[bit1] (b110) at (2.38,0) {1};
\node[bit0] (c110) at (2.86,0) {0};

\node[output1] (o110) at (2.38,-1.15) {1};

\draw[->,thick]
(2.38,-0.28) -- (o110.north);

\node[bit1] (a101) at (3.80,0) {1};
\node[bit0] (b101) at (4.28,0) {0};
\node[bit1] (c101) at (4.76,0) {1};

\node[output0] (o101) at (4.28,-1.15) {0};

\draw[->,thick]
(4.28,-0.28) -- (o101.north);

\node[bit1] (a100) at (5.70,0) {1};
\node[bit0] (b100) at (6.18,0) {0};
\node[bit0] (c100) at (6.66,0) {0};

\node[output0] (o100) at (6.18,-1.15) {0};

\draw[->,thick]
(6.18,-0.28) -- (o100.north);

\node[bit0] (a011) at (7.60,0) {0};
\node[bit1] (b011) at (8.08,0) {1};
\node[bit1] (c011) at (8.56,0) {1};

\node[output1] (o011) at (8.08,-1.15) {1};

\draw[->,thick]
(8.08,-0.28) -- (o011.north);

\node[bit0] (a010) at (9.50,0) {0};
\node[bit1] (b010) at (9.98,0) {1};
\node[bit0] (c010) at (10.46,0) {0};

\node[output0] (o010) at (9.98,-1.15) {0};

\draw[->,thick]
(9.98,-0.28) -- (o010.north);

\node[bit0] (a001) at (11.40,0) {0};
\node[bit0] (b001) at (11.88,0) {0};
\node[bit1] (c001) at (12.36,0) {1};

\node[output1] (o001) at (11.88,-1.15) {1};

\draw[->,thick]
(11.88,-0.28) -- (o001.north);

\node[bit0] (a000) at (13.30,0) {0};
\node[bit0] (b000) at (13.78,0) {0};
\node[bit0] (c000) at (14.26,0) {0};

\node[output0] (o000) at (13.78,-1.15) {0};

\draw[->,thick]
(13.78,-0.28) -- (o000.north);
\end{tikzpicture}

\caption{Illustration of the local rule for ECA 74.}
\label{fig:eca_rule74_local}
\end{figure}
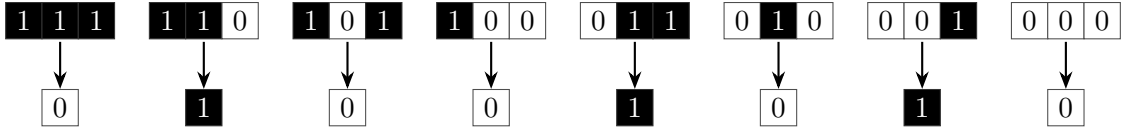

We consider the alphabet $\A:=\{0,1\}$ and the three consecutive integers $\{-1,0,1\}$
which define a discrete neighborhood $\NN_0$ of $0\in \Z$. The notation $\A^{\NN_0}$ represents the set
of the eight triples with entries in $\A$, and a \emph{local rule}
is a map $R\colon \A^{\NN_0}\to \A$ assigning to each of these triples an element of $\A$.
There exist exactly 256 local rules, which are called
\emph{Elementary Cellular Automata} (ECA). A typical illustration of one of these rules
is provided in Figure \ref{fig:eca_rule74_local}.
A systematic notation for these ECA has been introduced by Wolfram in the $80$'s, they are indexed by elements of the set \{0,1,\dots,255\}.
This notation allows one to deduce the local rule from the number assigned to each elementary cellular automaton. 

\begin{figure}[htbp]
   \centering
    \includegraphics[height=8cm]{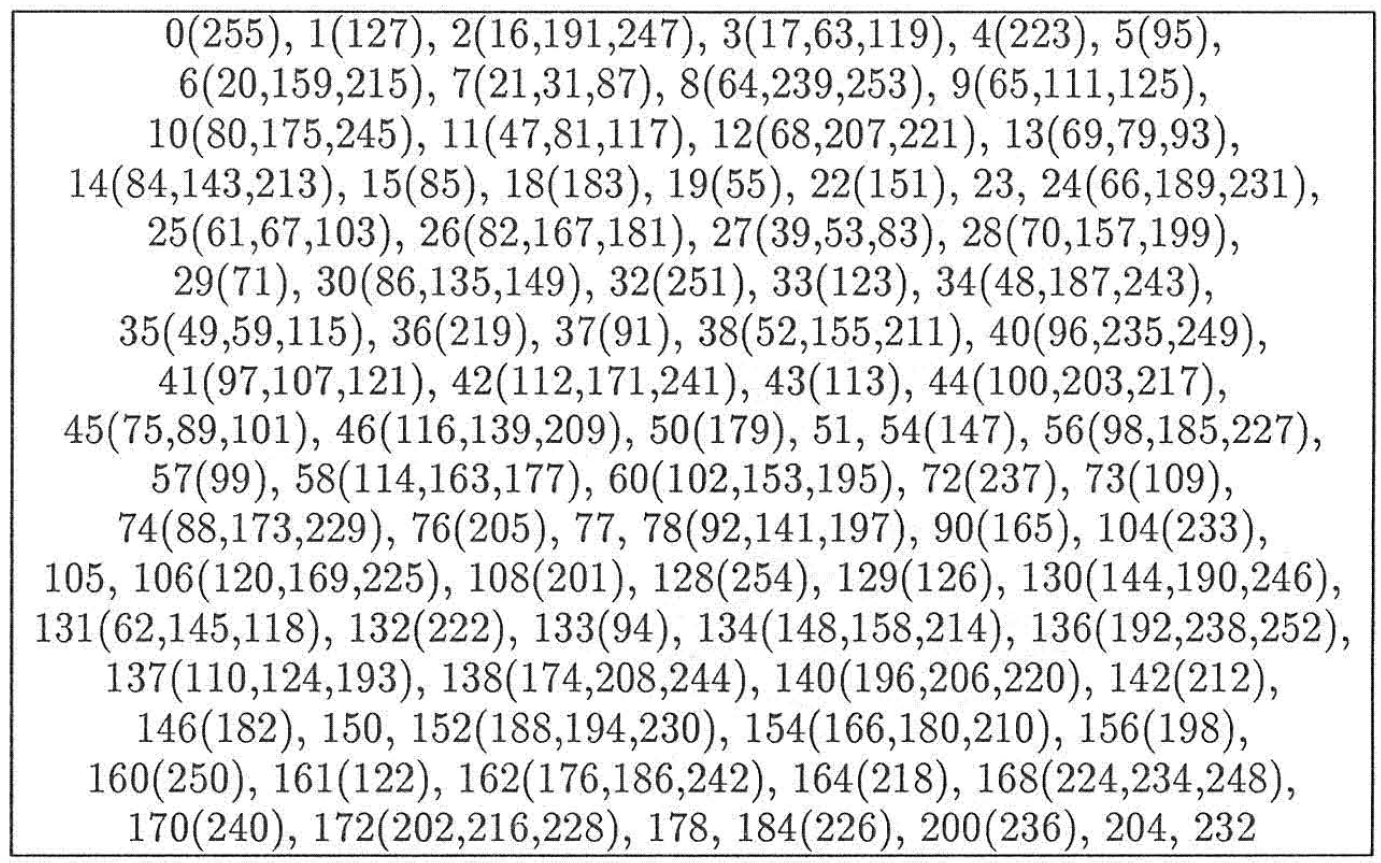}
    \caption{List of the 88 clusters of ECA, with the representative of each cluster chosen with the smaller decimal number,
    table reproduced from \cite{li1990structure}}
    \label{fig_LiP}
\end{figure}

By using two natural symmetries (complementary rule and mirrored rule), the 256 ECA can be classified into
88 clusters. The ECA inside each cluster are considered as equivalent, and the representative
of each cluster is chosen with the smaller decimal number, namely the one with the smaller index.
The list of these clusters is provided in Figure \ref{fig_LiP}. In the sequel, we investigate only the properties
of the representative of each cluster.

Each elementary cellular automaton can then be realized on a finite or an infinite set.
Since the infinite realization is out of the scope of this paper, we concentrate
on finite realizations. For that purpose, fix $n\in \N$ with $n\geq 4$, and consider
the $n$ integers $\NN:=\{1,2,\dots,n\}$. From the local rule, we define a \emph{global rule}, 
namely a map $R_n\colon \A^\NN\to \A^\NN$. In simpler words, the global rule associates to any $n$-tuples
with entries in $\A$ a new $n$-tuples with entries in $\A$. Clearly, the global rule should locally respect
the local rule. The precise formulation reads: for any $u\in \A^\NN$ and for 
$j\in \NN\setminus\{1,n\}$ we set
$$
[R_nu](j):=R\big(u(\cdot-j)|_{\NN_0}\big),
$$
which means that the value at $j$ is obtained by applying the local rule to the $n$-tuple shifted by $j$.
For $j\in \{1,n\}$, this procedure can not be applied, and a boundary condition has to be fixed.
Several solutions exist and this choice plays an important role, see \cite{LuValle} for an extensive analysis 
of the effect of boundary conditions on the resulting dynamical system.

For our investigations, we choose the most common boundary condition, namely the periodic boundary condition.
It corresponds to defining the global rule on a discrete circle made of $n$ points, and more precisely it means that
$$
[R_nu](1):=R\big(u(n), u(1), u(2)\big)\quad \hbox{ and }\quad  [R_nu](n):=R\big(u(n-1),u(n),u(1)\big).
$$
Our interest in the periodic boundary condition is that it does not involve any arbitrary choice, 
and corresponds to a natural realization on a discrete circle.
In the sequel, this choice is always tacitly assumed, and we call a \emph{finite periodic realization}
of the underlying ECA the map $R_n\colon \A^\NN\to \A^\NN$ endowed with periodic boundary condition.

For concreteness, for $n\in \N$ with $n\geq 4$, we call $u:=(u_1, u_2, \dots, u_n) \in \A^\NN$ a \emph{state}, and 
the set $\A^\NN$ the \emph{state space}. Clearly, the state space consists of $2^n$ states.
For any ECA, namely for any local rule $R$ indicated in Figure \ref{fig_LiP},
$R_n$ denotes the corresponding realization on the set $\NN$ with periodic boundary condition.
As a consequence, $R_n$ is a function from the state space to itself. 

\subsection{State transition network}\label{sec_stn}

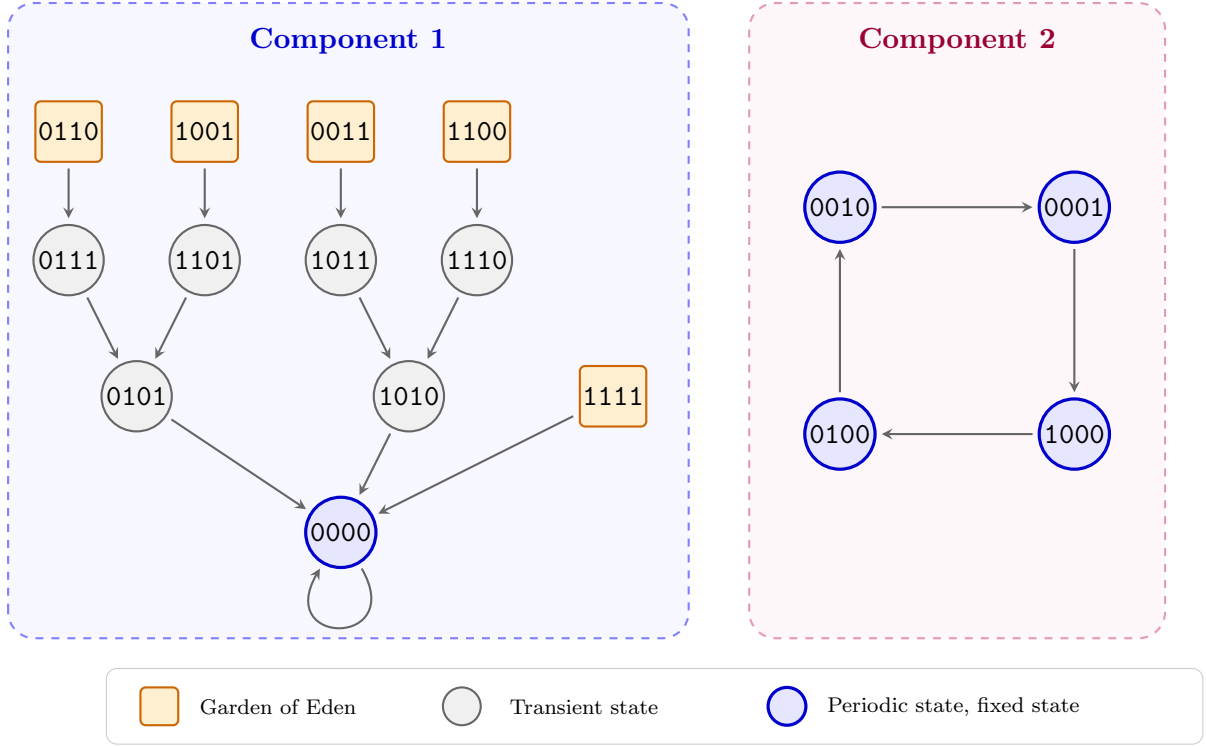
\begin{figure}
\begin{tikzpicture}[
    font=\small,
    >=stealth,
    base/.style={thick, draw=black, minimum size=8mm, inner sep=1pt, outer sep=2.5pt, font=\bfseries\ttfamily},
    fp_attractor/.style={base, circle, fill=blue!10, draw=blue!80!black, very thick},
    fp_transient/.style={base, circle, fill=gray!12, draw=black!60},
    fp_eden/.style={base, rectangle, fill=yellow!40!orange!20, draw=orange!80!black, rounded corners=2pt},
    cycle_attractor/.style={base, circle, fill=blue!10, draw=blue!80!black, very thick}]

\draw[dashed, draw=blue!50, fill=blue!3, rounded corners=10pt, thick] 
    (-0.8, -1.2) rectangle (8.2, 7.2);
\node[font=\bfseries, text=blue!80!black] at (3.7, 6.7) {Component 1};

\draw[dashed, draw=purple!40, fill=purple!3, rounded corners=10pt, thick] 
    (9.0, -1.2) rectangle (14.5, 7.2);
\node[font=\bfseries, text=purple!80!black] at (11.75, 6.7) {Component 2};

\node[fp_eden] (n0110) at (0.0, 5.5) {0110};
\node[fp_eden] (n1001) at (1.8, 5.5) {1001};
\node[fp_eden] (n0011) at (3.6, 5.5) {0011};
\node[fp_eden] (n1100) at (5.4, 5.5) {1100};
\node[fp_eden] (n1111) at (7.2, 2.0) {1111};

\node[fp_transient] (n0111) at (0.0, 3.8) {0111};
\node[fp_transient] (n1101) at (1.8, 3.8) {1101};
\node[fp_transient] (n1011) at (3.6, 3.8) {1011};
\node[fp_transient] (n1110) at (5.4, 3.8) {1110};

\node[fp_transient] (n0101) at (0.9, 2.0) {0101};
\node[fp_transient] (n1010) at (4.5, 2.0) {1010};

\node[fp_attractor] (n0000) at (3.6, 0.2) {0000};

\draw[->, thick, black!60] (n0110) -- (n0111);
\draw[->, thick, black!60] (n1001) -- (n1101);
\draw[->, thick, black!60] (n0011) -- (n1011);
\draw[->, thick, black!60] (n1100) -- (n1110);

\draw[->, thick, black!60] (n0111) -- (n0101);
\draw[->, thick, black!60] (n1101) -- (n0101);
\draw[->, thick, black!60] (n1011) -- (n1010);
\draw[->, thick, black!60] (n1110) -- (n1010);

\draw[->, thick, black!60] (n0101) -- (n0000);
\draw[->, thick, black!60] (n1010) -- (n0000);
\draw[->, thick, black!60] (n1111) -- (n0000);

\draw[->, thick, black!60] (n0000) to [loop below, in=-120, out=-60, distance=1.2cm] (n0000);

\node[cycle_attractor] (n0010) at (10.2, 4.5) {0010};
\node[cycle_attractor] (n0001) at (13.3, 4.5) {0001};
\node[cycle_attractor] (n1000) at (13.3, 1.5) {1000};
\node[cycle_attractor] (n0100) at (10.2, 1.5) {0100};

\draw[->, thick, black!60] (n0010) -- (n0001);
\draw[->, thick, black!60] (n0001) -- (n1000);
\draw[->, thick, black!60] (n1000) -- (n0100);
\draw[->, thick, black!60] (n0100) -- (n0010);

\begin{scope}[shift={(1.0, -2.2)}]

    \draw[draw=gray!40, fill=white, rounded corners=4pt] (-0.5, -0.4) rectangle (14, 0.6);
    
    \node[fp_eden, minimum size=5mm, font=\tiny] at (0.2, 0.1) {};
    \node[anchor=west, font=\scriptsize] at (0.6, 0.1) {Garden of Eden};
    
    \node[fp_transient, minimum size=5mm, font=\tiny] at (4.2, 0.1) {};
    \node[anchor=west, font=\scriptsize] at (4.7, 0.1) {Transient state};
    
    \node[fp_attractor, minimum size=5mm, font=\tiny] at (8.5, 0.1) {};
    \node[anchor=west, font=\scriptsize] at (8.9, 0.1) {Periodic state, fixed state};
\end{scope}
\end{tikzpicture}
    \caption{State transition network for ECA 74 and for $n=4$}
    \label{fig_74}
\end{figure}

Given a global rule $R_n$, the state space $\A^\NN$ is naturally endowed with the structure 
of a directed graph (with possibly several disjoint components).
The set of vertices of the graphs is given by the set of states, 
and edges are defined by the global rule, namely the map $u\mapsto R_n u$ defines a directed
edge from $u$ to $R_n u$ for any $u\in \A^\NN$. As a result, 
each vertex has an out-degree exactly equal to $1$. 
Thus, the graph does not contain multiple (directed) edges, but cycles and loops
certainly exist. The resulting directed graph is usually called \emph{state transition network}.
An example for the ECA $74$ and for $n=4$ is provided in Figure \ref{fig_74}.

Given this structure of directed graph, let us further decompose the set of vertices. 
Vertices with in-degree $0$ are called \emph{Garden of Eden}: they are vertices which 
are not the image by $R_n$ of any other vertex.
We call \emph{cyclic vertex} any vertex which belongs to a cycle (containing at least two vertices),
and a \emph{fixed vertex} a vertex which is the source (and the target) of a loop. Finally, any vertex which is neither
a Garden of Eden nor a cyclic or a fixed vertex is called a \emph{transient vertex}.
By extension, we use the same names for the corresponding states: Garden of Eden, cyclic state, fixed state, and transient state.

Before introducing a concrete example, let us mention an easy observation:
each component of a state transition network contains one and only one cycle or loop.
For later use, we denote by $C_m$ any cycle containing $m$ vertices. By extension, a fixed
vertex is denoted by $C_1$. As a consequence, to each component of the state
transition network is associated a unique $m\in \N$ corresponding
to the length of a unique cycle or fixed state.
The set 
\begin{equation}\label{eq_cycle_sp}
\GRn:=\big\{m\in \N\mid C_m  \hbox{ is contained in a component of the state transition network}\big\}
\end{equation}
is also called the \emph{cycle length spectrum} of the realization of the ECA $R$ on $\A^\NN$. 
Note that any value $m$ in $\GRn$ is repeated according to the appearance of a cycle $C_m$ (or fixed state $C_1$)
in the state transition network. As a consequence, the cardinality of $\GRn$ corresponds to the number
of components of the network.

\begin{example}\label{Ex_1}
We consider the ECA $74$ and fix $n=4$. In this case, the state space contains $16$ states
and the corresponding directed graph $16$ vertices. The state transition
network, represented in Figure \ref{fig_74}, consists of two components, 
one containing a fixed state and one consisting of a cycle $C_4$ only.
\end{example}

\subsection{Koopman operator}

We now introduce Koopman operator. In our setting,
this operator is going to be a linear operator acting on a Hilbert space.
Clearly, the Hilbert space depends on the parameter $n$ while Koopman
operator is going to depend both on the ECA $R$ and on $n$.

Let us consider all complex-valued functions defined on the phase space $\A^\NN$, 
namely all functions $f\colon \A^\NN\to \C$. Since the set $\A^\NN$ contains $2^n$
elements, these functions generate a vector space of dimension $2^n$.
Endowed with the standard scalar product $\langle \cdot,\cdot \rangle$ and $\ell^2$-norm $\|\cdot\|$, 
this set defines a finite dimensional Hilbert space which we denote by $\H_n$.
A standard basis of this Hilbert space is then given by the family
of Kronecker delta functions $\{\delta_u\}_{u\in \A^\NN}$.
Since $\H_n$ can be identified with $\C^{2^n}$, it follows 
that the set of all bounded operators acting on $\H_n$ can be identified with  
the set $M_{2^n}(\C)$ of $2^n\times 2^n$ complex matrices.

In the framework introduced above, the \emph{Koopman operator} associated with the realization
$R_n$ of the ECA $R$ is the operator defined by
\begin{equation}
\KRn\colon \H_n\to \H_n,
\qquad
[\KRn f](u) = f\big(R_n u\big) \quad \forall u \in \A^\NN.
\end{equation}
Clearly, Koopman operator defines a bounded and linear operator on $\H_n$,
which makes it of great interest. 
In the standard basis of $\H_n$, this operator can be represented by an
element of $M_{2^n}(\C)$ by computing its matrix elements, namely
\begin{align*}
\langle \delta_v,\KRn \delta_u\rangle 
& = \sum_{w\in \A^\NN} \overline{\delta_v(w)} [\KRn \delta_u](w) \\
& = [\KRn \delta_u](v) \\
&= \delta_u\big(R_n v\big) \\
&= \begin{cases}1 & \hbox{ if }R_n v=u , \\ 0 & \hbox{ otherwise.}\end{cases}
\end{align*}

Having defined Koopman operator, one can look at its set of eigenvalues.
As emphasized in the Introduction, our interest is to study the ECA 
according to the spectra $\sigma(\KRn)$ of Koopman operators $\KRn$ for various $n$, since this approach
does not involve the choice of any specific initial conditions. 
Thus, one might be reduced to studying the finite-dimensional
eigenvalue equation
\begin{equation*}
\KRn f = \lambda f
\end{equation*}
for $f \in \H_n$ with $f\neq 0$ and for $\lambda \in \C$.
However, this approach would quickly be impracticable since the size
of the matrices would become too large for large $n$. 
Fortunately, a link between the spectrum of $\KRn$ and the state-transition
network has been established in \cite{Taga1}, and this link
makes the computations much more efficient.
In the next statement, we summarize the content of Lemmas 1 to 6
of \cite{Taga1}. Recall that the cycle length spectrum $\GRn$
has been introduced in \eqref{eq_cycle_sp}.

\begin{prop}\label{prop_spectrum}
1) For any $m\in \GRn$, the operator $\KRn$ possesses the eigenvalues $\e^{2\pi i k/m}$
for any $k \in \{1, \dots,m\}$,

2) To any state which is a Garden of Eden or a transient state corresponds an eigenvalue
$0$ of the operator $\KRn$.
\end{prop}

Observe that the content of this statement provides the $2^n$ eigenvalues of the operator
$\KRn$. Indeed, each cycle $C_m$ (or fixed state for $m=1$) produces $m$ eigenvalues, and each
transient state or Garden of Eden corresponds to one eigenvalue. By Summing these
contributions, one gets $2^n$ eigenvalues, as expected. One also deduces that the multiplicity
of the eigenvalue $1$ is equal to the number of components of the state-transition network,
see \cite[Thm.~2]{Taga1}. Similarly, the multiplicity of $0$ coincide with the number
of states which do not belong to any periodic orbit or which are not fixed states,
as emphasized in \cite[Thm.~1]{Taga1}. 
Let us finally stress that the spectrum of $\KRn$ consists only of two types of eigenvalues:
either they are are of modulus $1$, or are $0$. This is a direct outcome of 
Proposition \ref{prop_spectrum}.

\begin{example}\label{Ex_2}
Coming back to Example \ref{Ex_1} for ECA $74$ and $n=4$, the spectrum of 
$K_{74_4}$ can be deduced from Figure \ref{fig_74}, once Proposition \ref{prop_spectrum}
is taken into account. Indeed, one directly gets:
\begin{equation}
\sigma\big(K_{74_4}\big)=\Big\{
\underbrace{1}_{\text{due to }C_1}, \underbrace{i,-1,-i,1}_{\text{due to } C_4},
\underbrace{0,0,0,0,0,0,0,0,0,0,0}_{\text{due to Garden of Eden and transient states}}\Big\}.
\end{equation}
\end{example}

Let us close this section with a remark: Koopman operator contains much more information
than its mere spectrum. For example, its Jordan block decomposition is directly related to 
the graph structure exhibited in Section \ref{sec_stn}. However, as already emphasized these matrices
are quickly becoming too big, as $n$ increases. On the other hand, getting the information
about the cycle length spectrum $\GRn$ is relatively easy, and therefore computing the spectrum
of Koopman operator does not take too much resources.

\section{Classification}\label{sec_classification}

Based on the spectra of Koopman operators, we establish a classification of all ECA.
Let us immediately stress that this classification is based only on the characteristic of the spectra,
namely one the set $\sigma(\KRn)$ of eigenvalues of $\KRn$ including their multiplicities.
This information can be deduced from the cycle length spectrum, as already emphasized.
Note that the eigenvalue $0$ does not play a major role for the classification, 
but its minor role will be mentioned subsequently.
Note also that the full structure of the directed graph introduced in Section \ref{sec_stn}
is not taken into account, since this structure becomes too complicated for large $n$.
Clearly, this classification is partially heuristic, but nevertheless it provides useful information.

\begin{remark}
Since our classification is based on the spectra of Koopman operators, it is important to 
represent several of them. As shown in Proposition \ref{prop_spectrum}, the spectrum consists only of eigenvalues
on the unit circle of the complex plane, and of the eigenvalue $0$. In the subsequent figures,
we represent the non-zero eigenvalues together with their multiplicity. For that purpose,
observe that any $\lambda \in \C$ with $|\lambda|=1$ can be written uniquely as
$\lambda = \e^{i\theta}$ for some $\theta \in (0,2\pi]$. We call
$\theta$ the argument of $\lambda$, and have chosen the cut at $\theta=0$ in order
to have the eigenvalue $\lambda=1$ represented at $\theta = 2\pi$.
Obviously, the eigenvalue $-1$ is represented at $\theta=\pi$, 
while no eigenvalue is represented at $\theta=0$ since since value is excluded.

On the same pictures we also give an information about the $0$-eigenvalue by
providing the \emph{zero eigenvalue ratio}, namely the ratio of the multiplicity of 
the eigenvalue $0$ divided by the total number of eigenvalues $2^n$, see \eqref{eq_ratio}
and the explanations provided there.
\end{remark}

The 88 clusters of ECA have been divided into five main categories, 
with some natural subdivisions inside some of them.
For concreteness, our analysis is motivated by the realizations for $n\in \{13, 14, \dots, 24\}$, 
with the spectrum computed from $\GRn$ as explained in Proposition \ref{prop_spectrum}.
The information on $\GRn$, namely the list of all cycles $C_m$ with their multiplicity is provided in \cite{Zhu}. 
For each fixed $n$, let us also denote by $\Cmax(n)$ the largest value
of $\GRn$, \emph{i.e.}~the length of the largest cycle of the transition state network.
For shortness, we shall also speak about the \emph{spectra of the ECA $R$} while the precise meaning
is \emph{the spectra of Koopman operators $\KRn$ for the ECA $R$}.

\begin{figure}[htbp]
\centering
\begin{minipage}[b]{0.5\textwidth}
\centering
\includegraphics[width=0.9\textwidth]{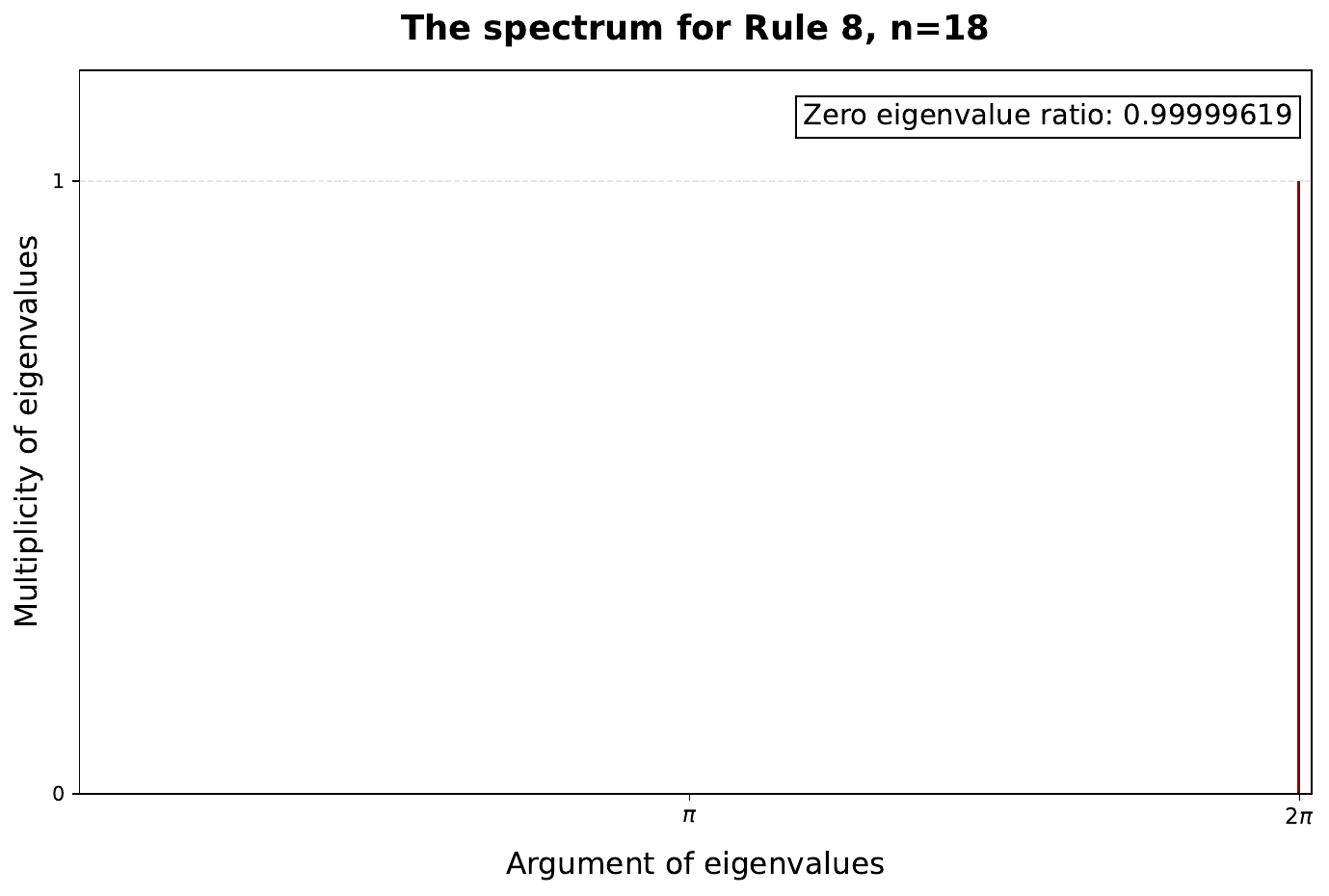}
\end{minipage}%
\begin{minipage}[b]{0.5\textwidth}
\centering
\includegraphics[width=0.9\textwidth]{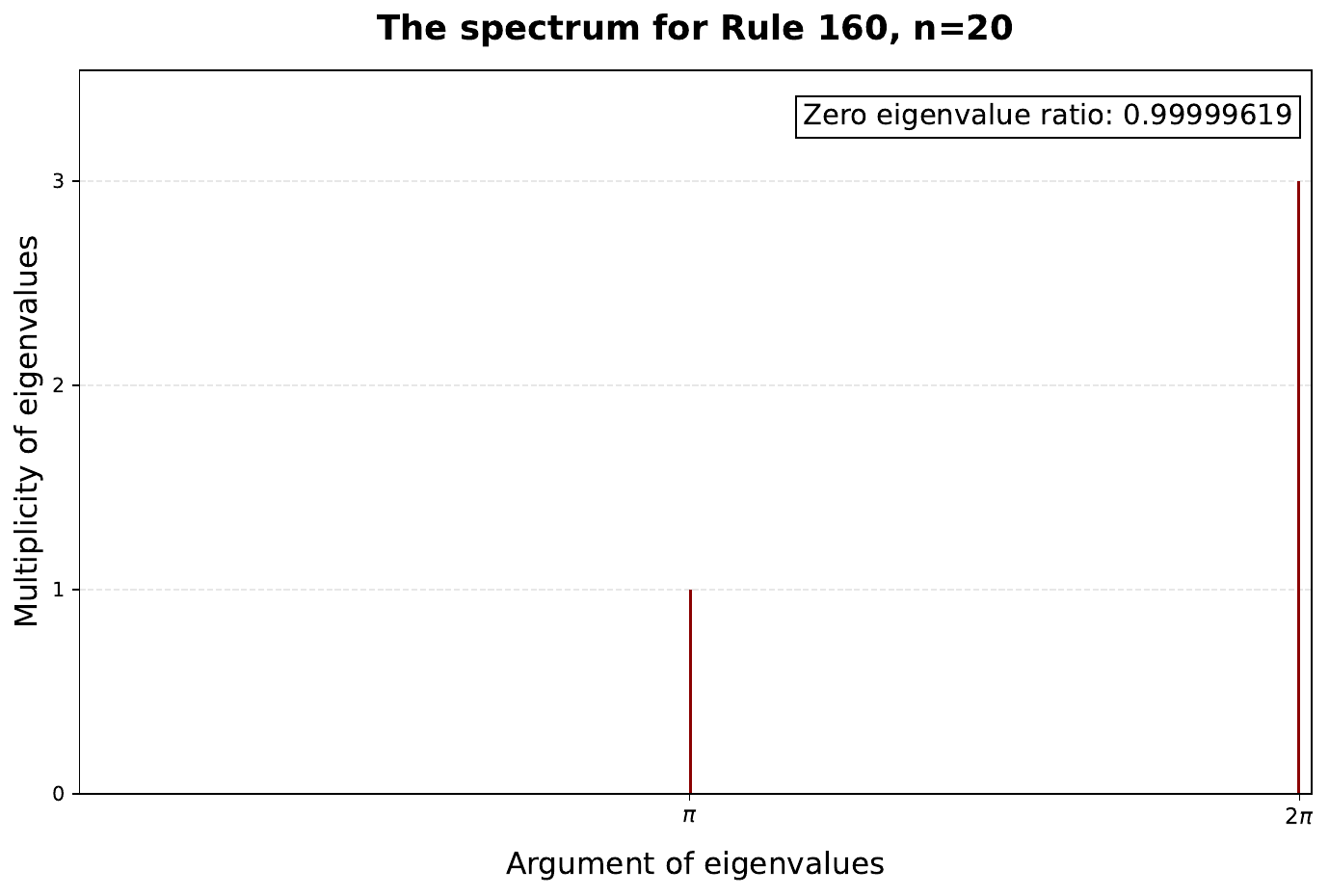}
\end{minipage}%
\caption{Two typical spectra of Category I: Rule $8$ for $n=18$ and Rule $160$ for $n=20$}
\label{fig_Cat_I}
\end{figure}

\paragraph{Category I}
This category contains all ECA whose spectra (excluding $0$) are included in $\{-1,1\}$ and with small multiplicity,
see Figure \ref{fig_Cat_I}.
For them, the cycle length spectrum $\GRn$ consists only of $1$ and $2$, each of them up to multiplicity $2$. 
The corresponding ECA are
$$
0, 8, 32, 128, 136, 160.
$$

Category I contains naturally three families: 
\begin{enumerate}
\item[I${}_1$] corresponding to ECA whose spectra (excluding $0$) is composed only of $1$ with multiplicity $1$
due to a single fixed vertex in the state transition network: Rules $0$ and $8$, see Figure \ref{fig_Cat_I} (left),
\item[I${}_2$] corresponding to ECA whose spectra (excluding $0$) is composed only of $1$ with multiplicity $2$ 
due to two fixed vertices in the state transition network: Rules $128$ and $136$,
\item[I${}_{1/2}$] corresponding to ECA whose spectra (excluding $0$) is composed only of $1$ of multiplicity $1$ for $n$ odd, 
and of $\{-1,1\}$ of multiplicity at most $3$ for $n$ even: Rules $32$ and $160$, see Figure \ref{fig_Cat_I} (right).
\end{enumerate}

\paragraph{Category II}
This category contains all ECA whose spectra contains $-1$ or $1$ (or both) with large multiplicities.
A few other additional eigenvalues (excluding $0$) might also be present, but only with small multiplicity
(here the ratio between smaller and larger values is at least of a factor $10$). 
An additional ECA also belongs to this category: Rule $62$, which has the $3$-rd complex root of $1$ in its spectrum 
with large multiplicities, and a few additional eigenvalues (excluding $0$) with small multiplicity.
The corresponding ECA are
\begin{align*}
& 1, 4, 5, 12, 13, 19, 23, 28, 29, 33, 36, 37, 44, 50, 51, 62, 72, 76, 77, 78, \\
& 104, 108, 132, 140, 156, 164, 172, 178, 200, 204, 232.
\end{align*}

\begin{figure}[htbp]
\centering
\begin{minipage}[b]{0.5\textwidth}
\centering
\includegraphics[width=0.9\textwidth]{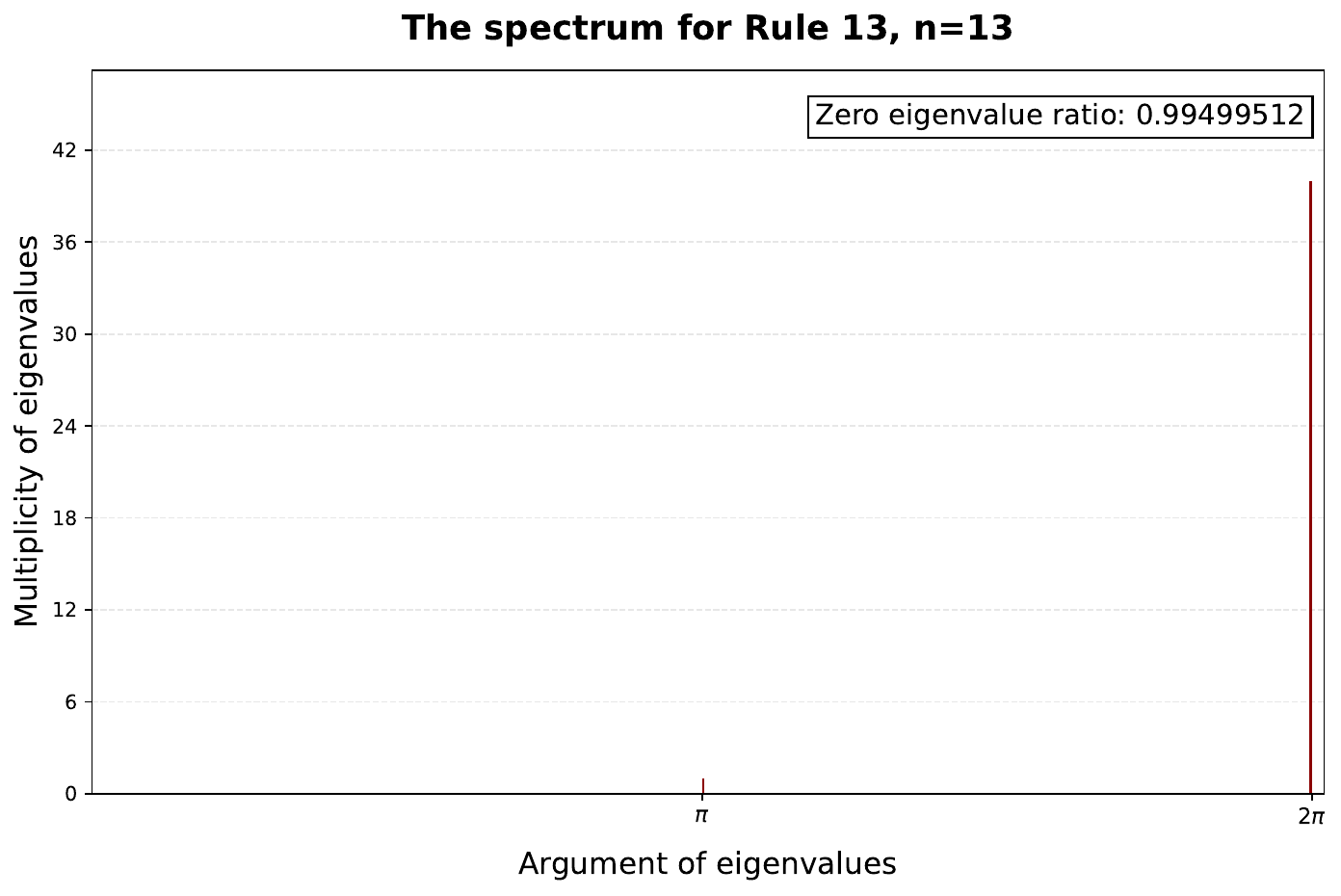}
\label{fig_Rule13_n13}
\end{minipage}%
\begin{minipage}[b]{0.5\textwidth}
\centering
\includegraphics[width=0.9\textwidth]{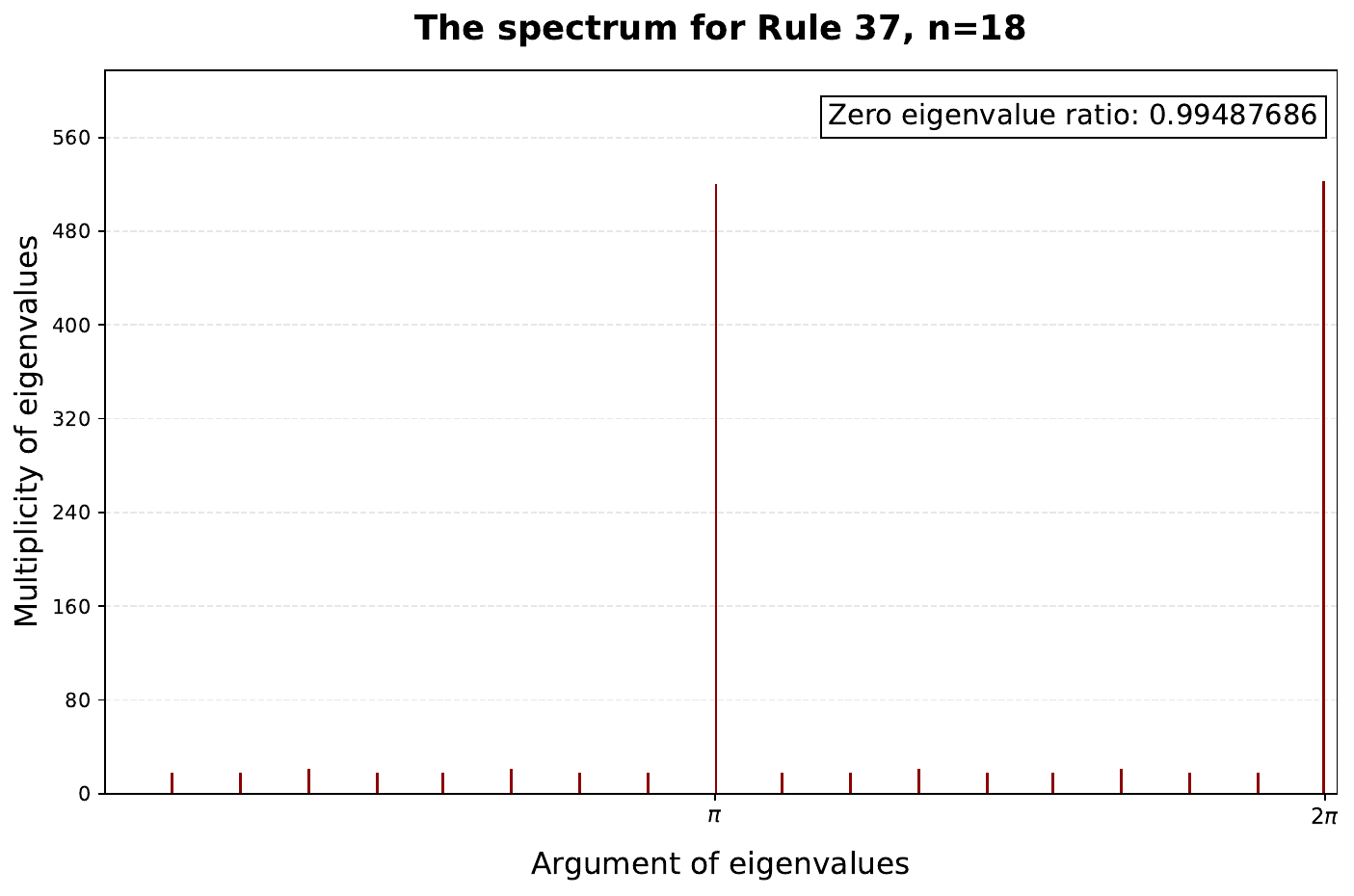}
\label{fig_Rule37_n18}
\end{minipage}%
\caption{Two typical spectra of Category II: Rule $13$ for $n=13$ and Rule $37$ for $n=18$}
\label{fig_Cat_II_1}
\end{figure}

\begin{figure}[htbp]
\centering
\begin{minipage}[b]{0.5\textwidth}
\centering
\includegraphics[width=0.9\textwidth]{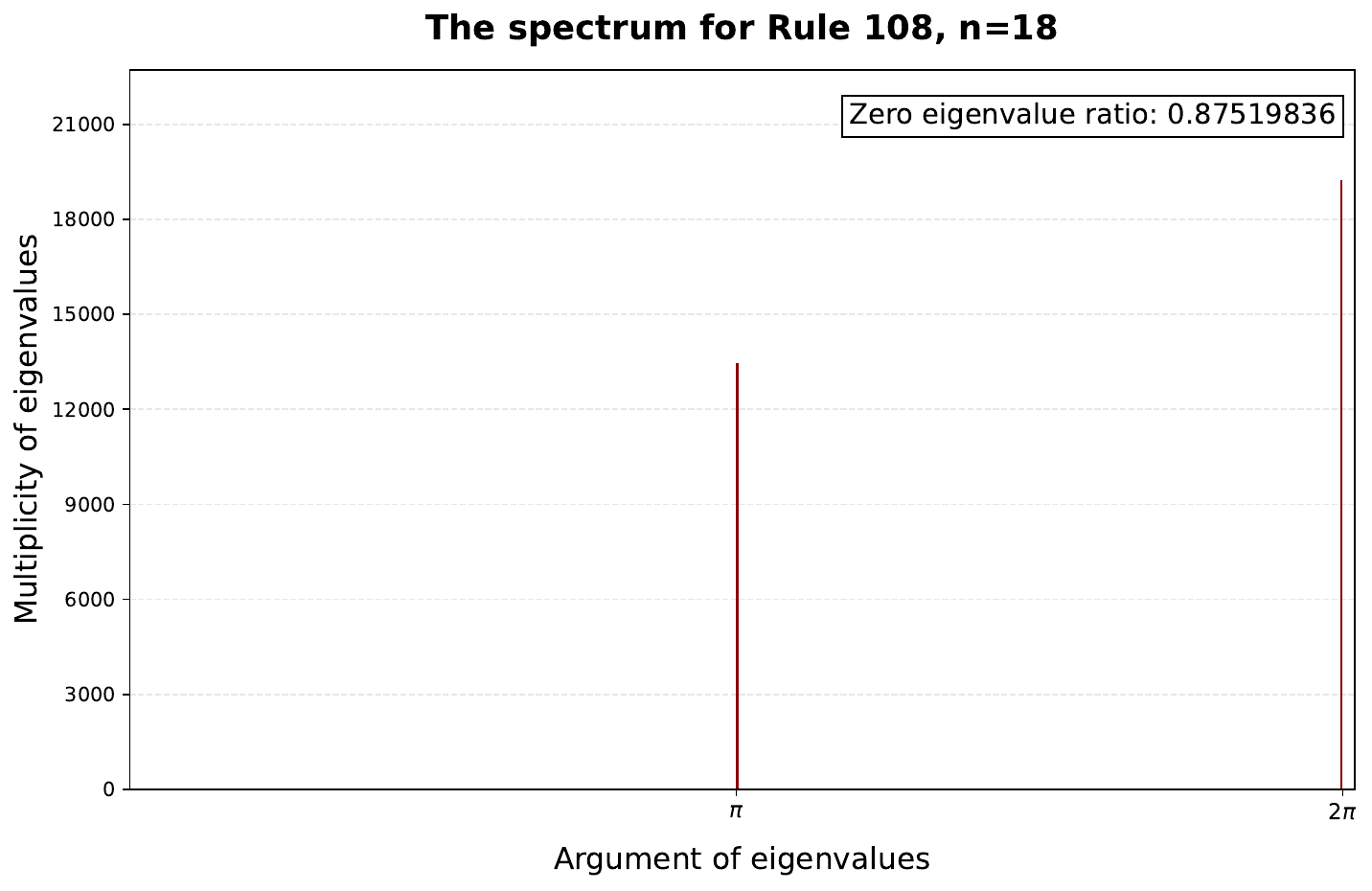}
\end{minipage}%
\begin{minipage}[b]{0.5\textwidth}
\centering
\includegraphics[width=0.9\textwidth]{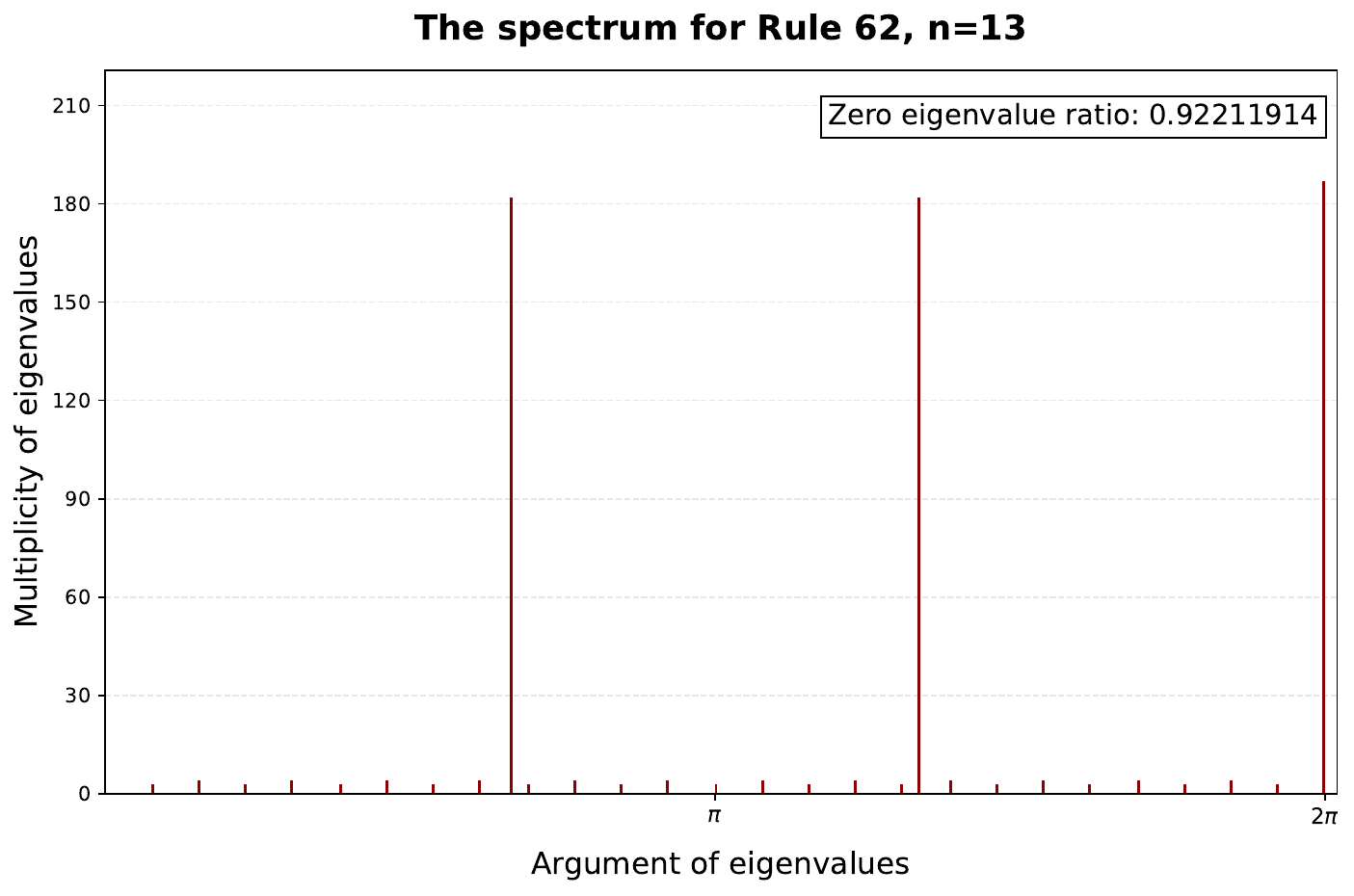}
\end{minipage}%
\caption{Two exceptional spectra of Category II: Rule $108$ for $n=18$ and Rule $62$ for $n=13$}
\label{fig_Cat_II_2}
\end{figure}

Category II contains naturally two families and two exceptional ECA. 
\begin{enumerate}
\item[II${}_1$] corresponding to ECA whose spectra (excluding $0$) is composed of $1$ with large multiplicity, 
and possibly of a few additional eigenvalues of small multiplicity: Rules $4$, $12$, $13$, $36$, $44$,  $72$, $76$, 
$77$, $78$, $104$, $132$, $140$, $164$, $172$, $200$, $204$, $232$, see Figure \ref{fig_Cat_II_1} (left),
\item[II${}_2$] corresponding to ECA whose spectra (excluding $0$) is composed of $-1$ and $1$ with similar multiplicities 
due to $2$-periodic orbits in the state transition network. A few additional eigenvalues of small multiplicity can also appear:  
Rules $1$, $5$, $19$, $23$, $28$, $29$, $33$, $37$, $50$, $51$, $156$, $178$, see Figure \ref{fig_Cat_II_1} (right),
\item[II${}_{1/2}$] containing an exceptional ECA whose state transition network contains several fixed vertices 
and several $2$-periodic orbits, generating a spectrum (excluding $0$) made of $-1$ and $1$ with large multiplicity: 
Rule $108$, see Figure \ref{fig_Cat_II_2} (left),
\item[II${}_3$] containing the exceptional ECA having the $3$-rd complex root of $1$ in its spectrum with large multiplicity, 
and a few additional eigenvalues (excluding $0$) with small multiplicity. Rule $62$, see Figure \ref{fig_Cat_II_2} (right).
\end{enumerate}

\paragraph{Category III}
This category contains all ECA whose cycle length spectrum $\GRn$ contains a large number
of maximal cycles of size $n$ or $2n$. Accordingly, the spectrum consists
mainly of the $n$-th complex roots of unity, or the $2n$-th complex roots of unity.
The cycle length spectrum $\GRn$ might also contain a small number of additional cycles of length corresponding to divisors of $n$ or of $2n$.
The corresponding ECA are 
\begin{align*}
& 2, 3, 6, 7, 10, 11, 14, 15, 24, 27, 34, 35, 38, 40, 42, 43, 46, 56, 57, 58, \\ 
& 130, 134, 138, 142, 152, 162, 168, 170, 184. 
\end{align*}

\begin{figure}[htbp]
\centering
\begin{minipage}[b]{0.5\textwidth}
\centering
\includegraphics[width=0.9\textwidth]{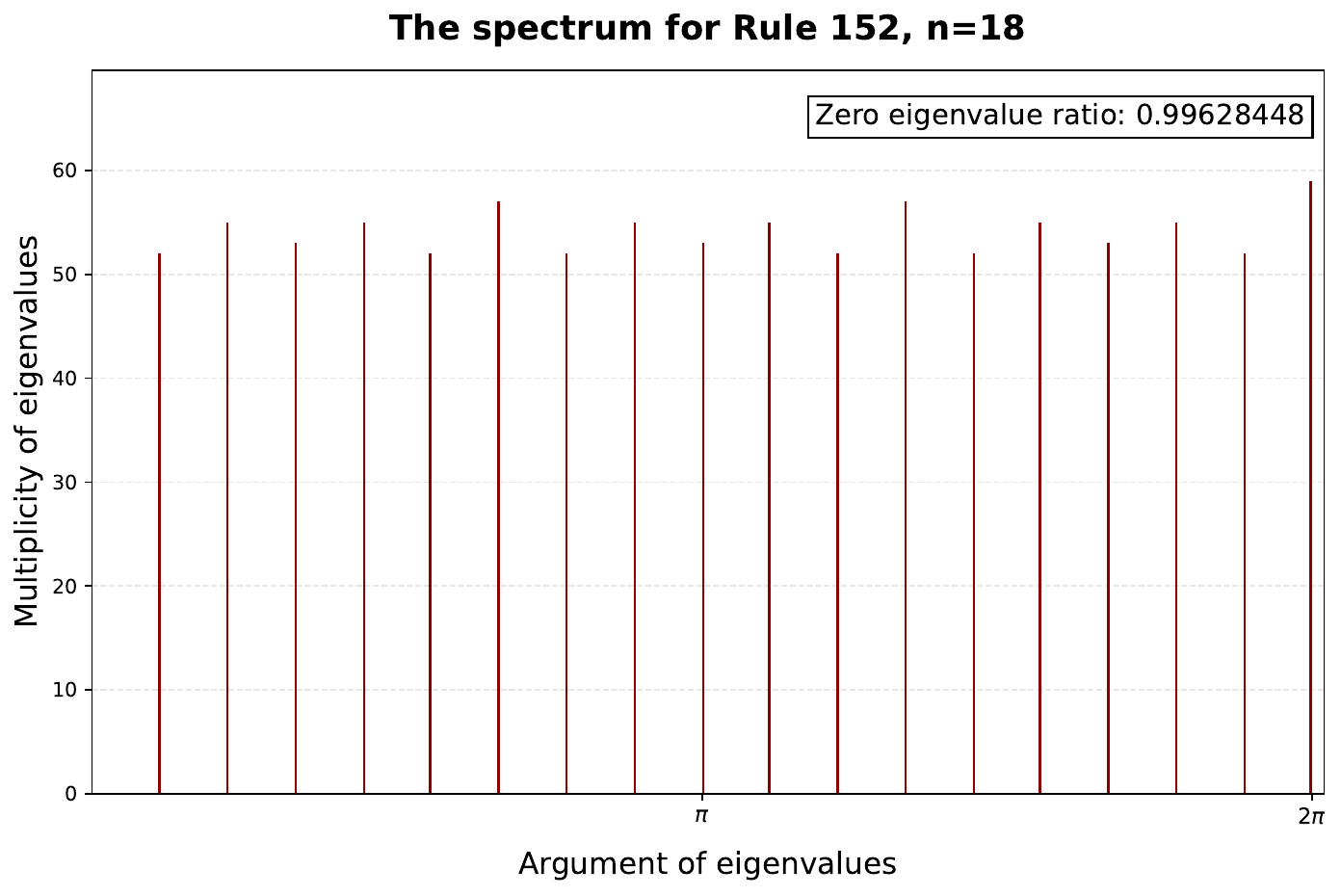}
\end{minipage}%
\begin{minipage}[b]{0.5\textwidth}
\centering
\includegraphics[width=0.9\textwidth]{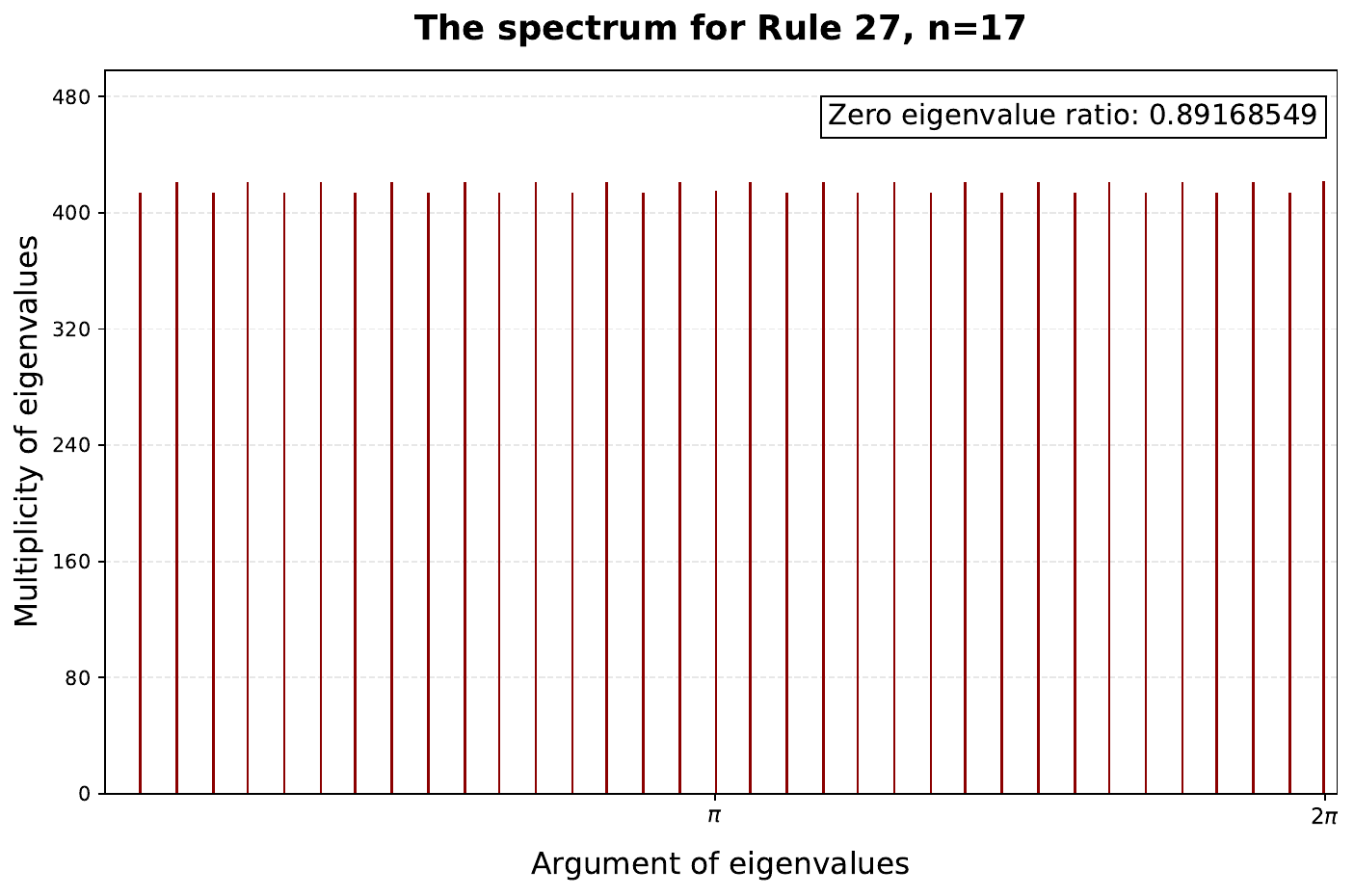}
\end{minipage}%
\caption{Two typical spectra of Category III: Rule $152$ for $n=18$ and Rule $27$ for $n=17$}
\label{fig_Cat_III_1}
\end{figure}

\begin{figure}[htbp]
\centering
\begin{minipage}[b]{0.5\textwidth}
\centering
\includegraphics[width=0.9\textwidth]{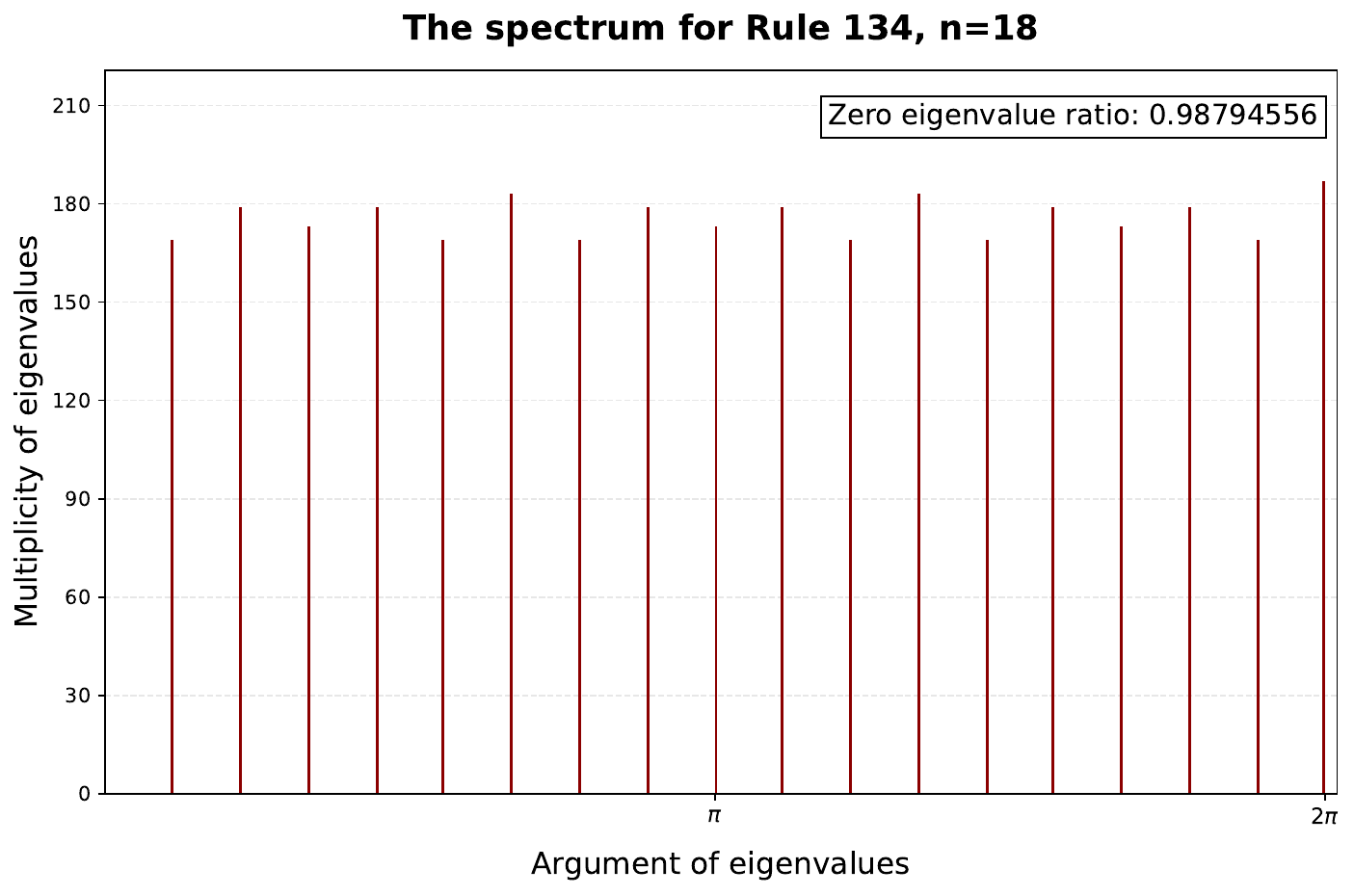}
\end{minipage}%
\begin{minipage}[b]{0.5\textwidth}
\centering
\includegraphics[width=0.9\textwidth]{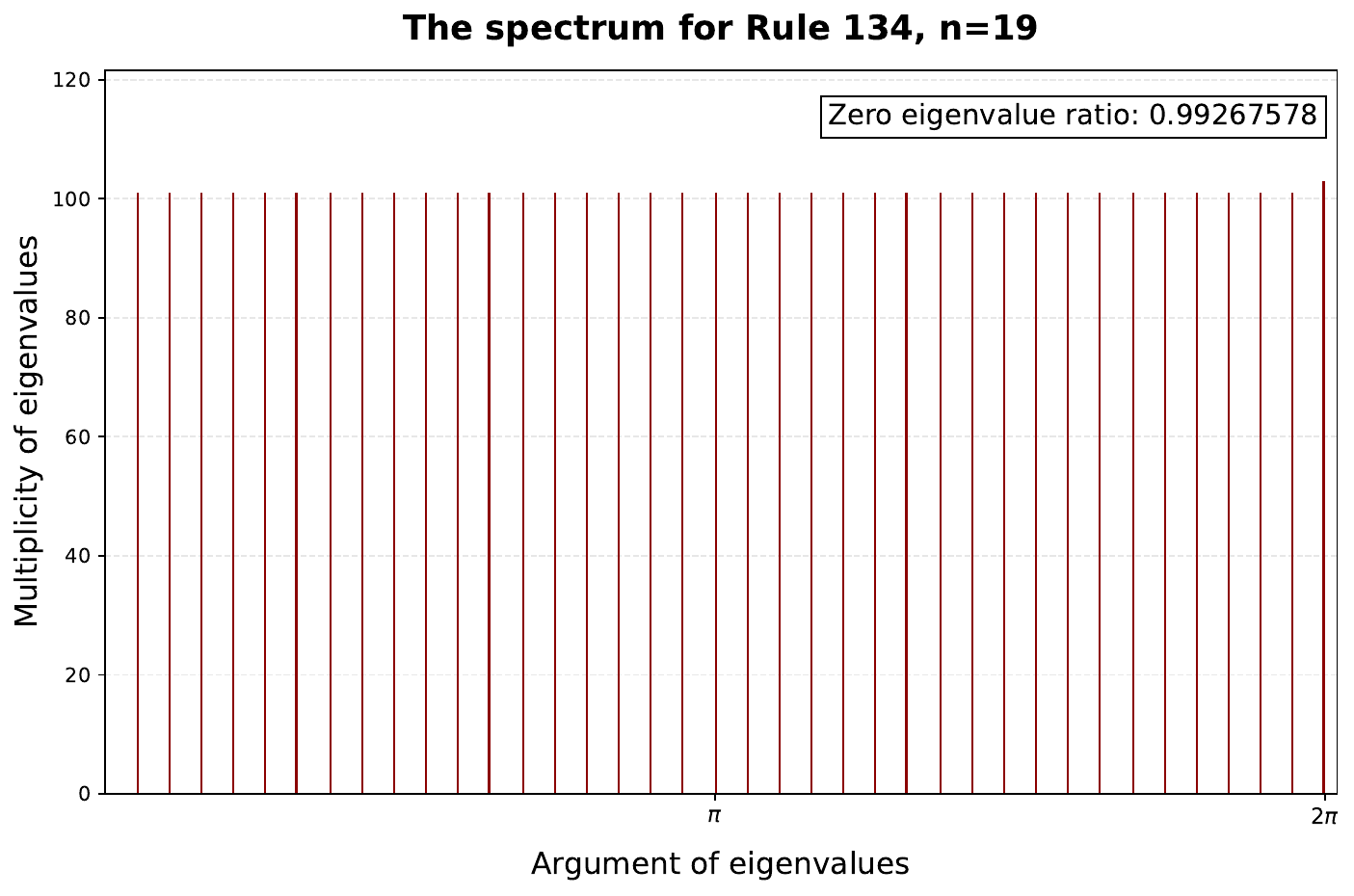}
\end{minipage}%
\caption{Two typical spectra of Category III${}_{1/2}$: Rule $134$ for $n=18$ and $n=19$}
\label{fig_Cat_III_2}
\end{figure}

Category III contains naturally three families of ECA. 
\begin{enumerate}
\item[III${}_1$] corresponding to ECA whose cycle length spectrum $\GRn$ contains
several maximal cycles of size $n$, leading to a spectrum consisting of $n$-th complex roots of unity with large multiplicity.  
A few additional eigenvalues of small multiplicity can also appear
and correspond to cycles in the state transition network of length dividing $n$:
Rules $2$, $10$, $24$, $34$, $40$, $42$, $46$, $56$, $57$, $130$, $138$, $152$, $162$, 
$168$, $170$, $184$, see Figure \ref{fig_Cat_III_1} (left),
\item[III${}_2$] corresponding to ECA whose cycle length spectrum $\GRn$ contains
several maximal cycles of size $2n$, leading to a spectrum consisting of $2n$-th complex roots of unity with large multiplicity.  
A few additional eigenvalues of small multiplicity can also appear
and correspond to cycles in the state transition network of length dividing $2n$:
Rules $3$, $7$, $27$, $35$, $58$, see Figure \ref{fig_Cat_III_1} (right),
\item[III${}_{1/2}$] corresponding to ECA whose cycle length spectrum $\GRn$ contains
many maximal cycles of size $n$ for $n$ even, and many maximal cycles of size $2n$ for $n$ odd. 
This leads to a spectrum consisting of $n$-th complex roots of unity with large multiplicity for $n$ even, 
and to a spectrum consisting of $2n$-th complex roots of unity with large multiplicity for $n$ odd. 
A few additional eigenvalues of small multiplicity can also appear
and correspond to cycles in the state transition network of length dividing $n$ or $2n$ respectively:
Rules $6$, $11$, $14$, $15$, $38$, $43$, $134$, $142$, see Figure \ref{fig_Cat_III_2} (left and right).
\end{enumerate}

\paragraph{Category IV}
This category contains all ECA whose cycle length spectrum $\GRn$ is rather
unpredictable, when looking at different $n$. 
It contains large numbers, but not extremely large. 
The corresponding ECA are
$$
9, 18, 22, 25, 26, 41, 54, 60, 73, 74, 90, 94, 105, 110, 122, 126, 146, 150, 154.
$$

\begin{figure}[htbp]
\centering
\begin{minipage}[b]{0.5\textwidth}
\centering
\includegraphics[width=0.9\textwidth]{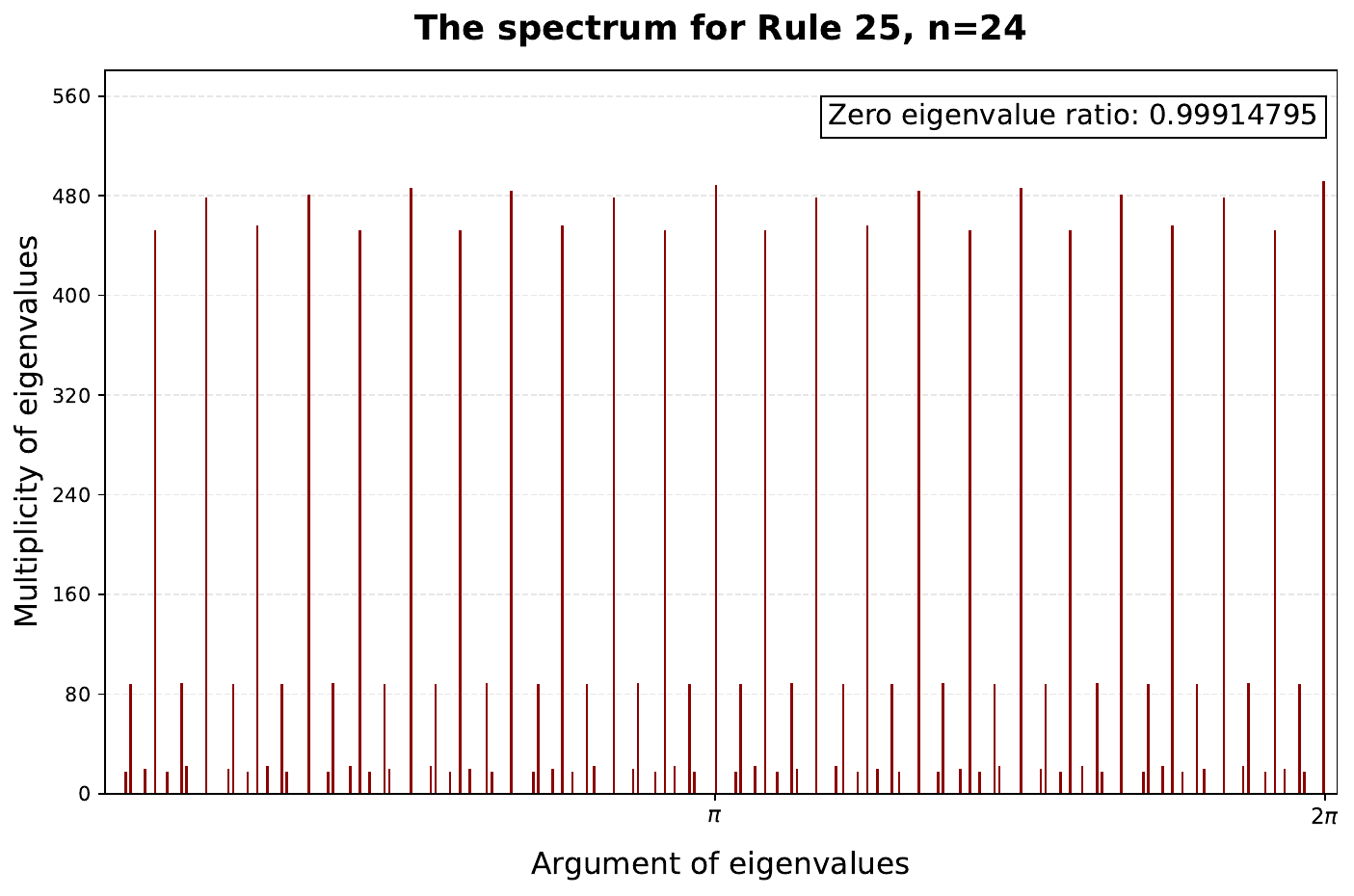}
\end{minipage}%
\begin{minipage}[b]{0.5\textwidth}
\centering
\includegraphics[width=0.9\textwidth]{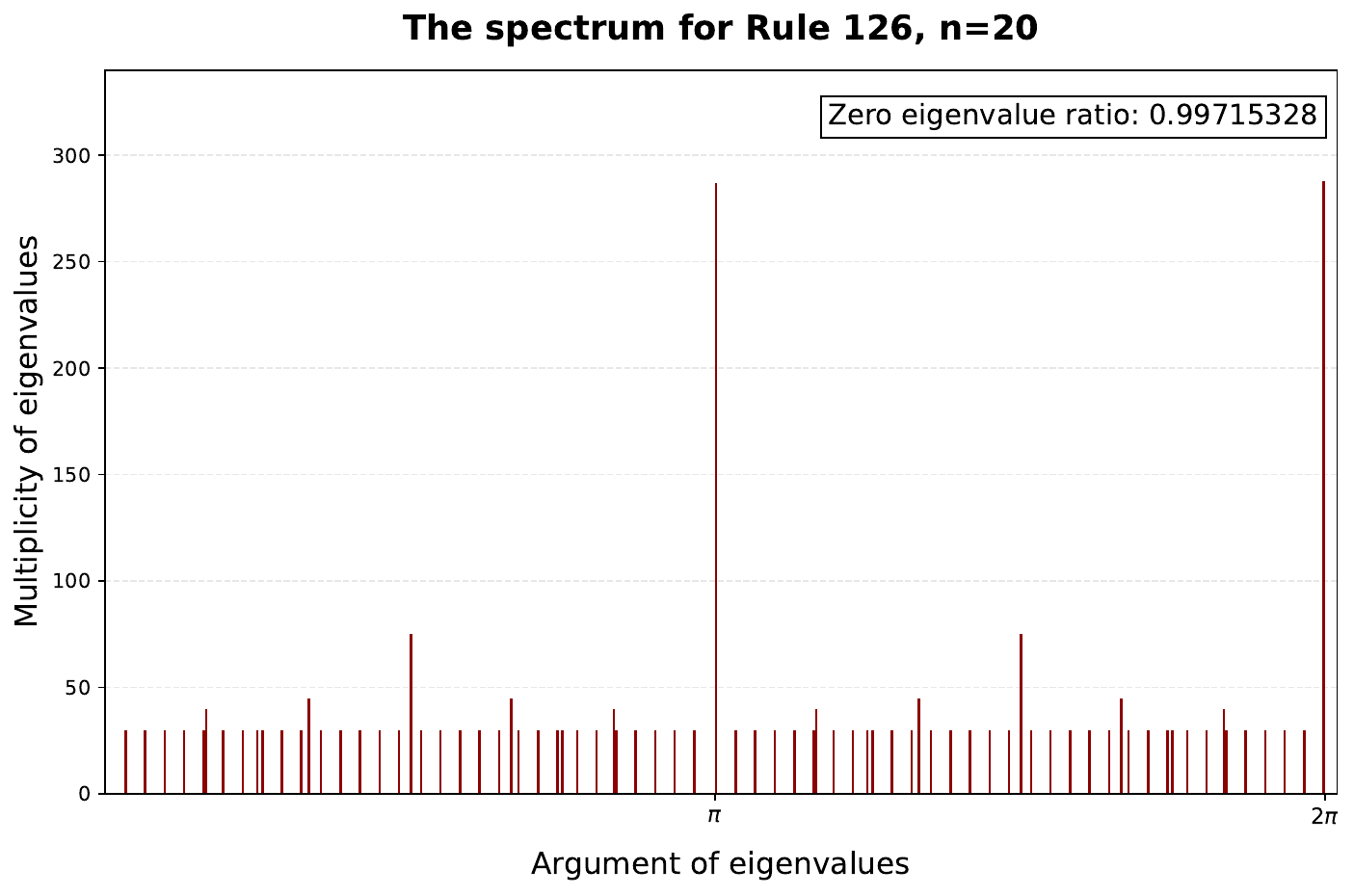}
\end{minipage}%
\caption{Two typical spectra of Category IV: Rule $25$ for $n=24$ and Rule $126$ for $n=20$}
\label{fig_Cat_IV_1}
\end{figure}

Two examples of spectra for Category IV are presented in Figure \ref{fig_Cat_IV_1}.
However, dividing Category IV into families is more delicate, and requires additional ideas introduced in the next section.
We shall then come back to this category.

\begin{figure}[htbp]
\centering
\begin{minipage}[b]{0.5\textwidth}
\centering
\includegraphics[width=0.9\textwidth]{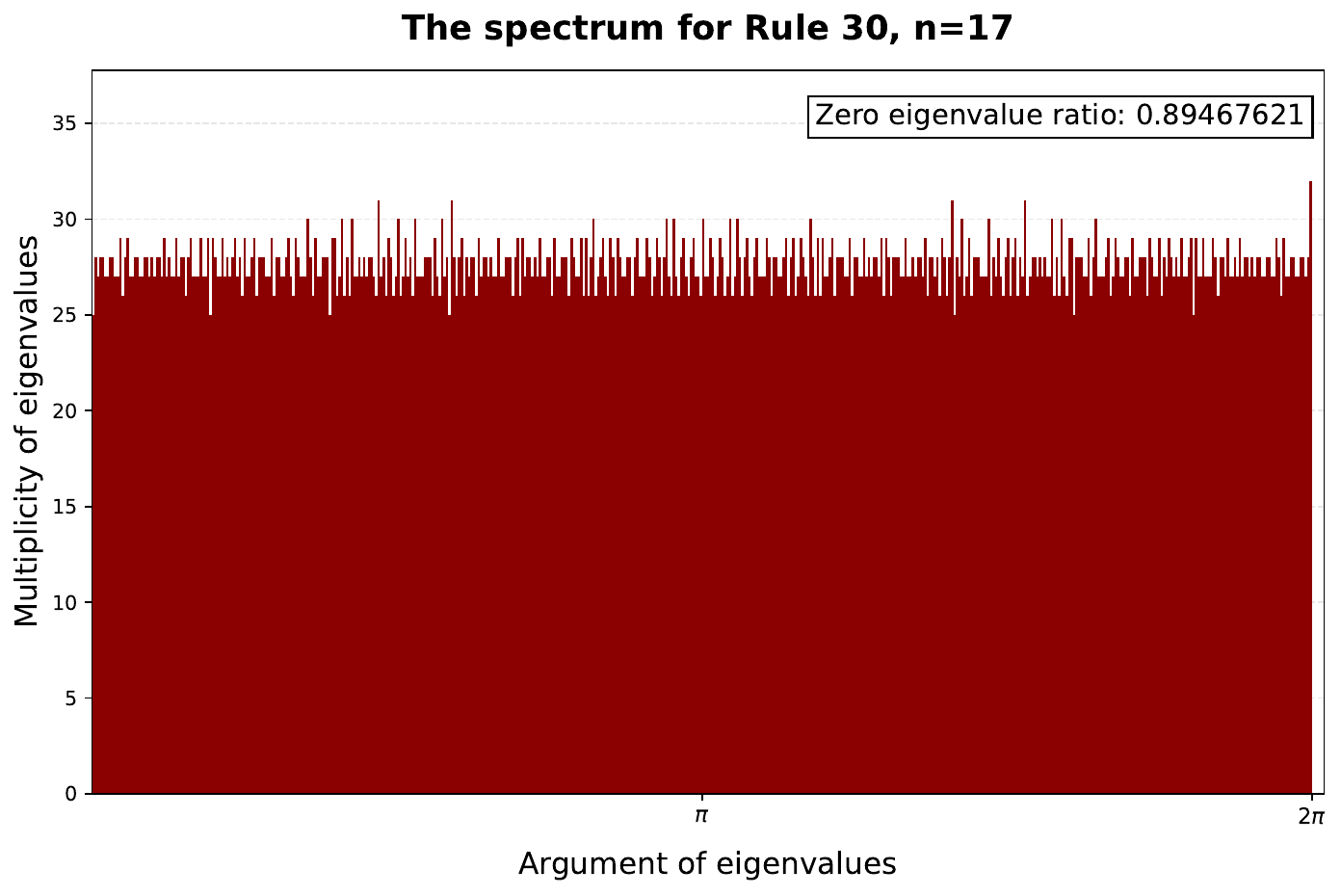}
\end{minipage}%
\begin{minipage}[b]{0.5\textwidth}
\centering
\includegraphics[width=0.9\textwidth]{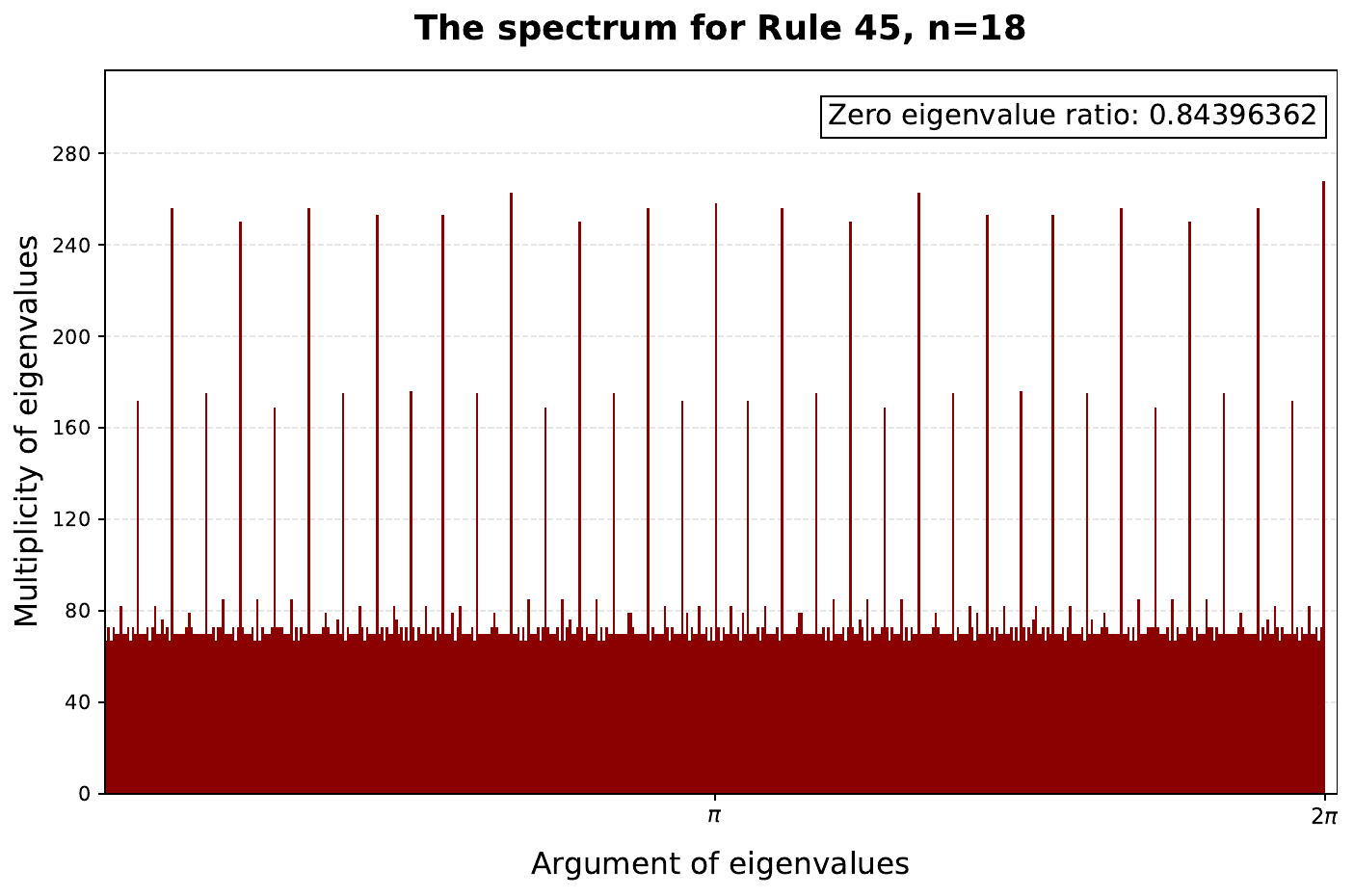}
\end{minipage}%
\caption{Two typical spectra of Category V: Rule $30$ for $n=17$ and Rule $45$ for $n=18$}
\label{fig_Cat_V}
\end{figure}

\paragraph{Category V}
This category contains all ECA whose cycle length spectrum $\GRn$ contains extremely large
numbers, or accordingly some extremely large $\Cmax(n)$. These cycles are usually of multiplicity $1$ or $2$. As a consequence, the spectrum of $\KRn$ looks like a dense set 
in the complex unit circle. The corresponding ECA are
$$
30, 45, 106. 
$$
Two examples of spectra for Category V are presented in Figure \ref{fig_Cat_V}.

\begin{figure}[htbp]
   \centering
    \includegraphics[height=5.2cm]{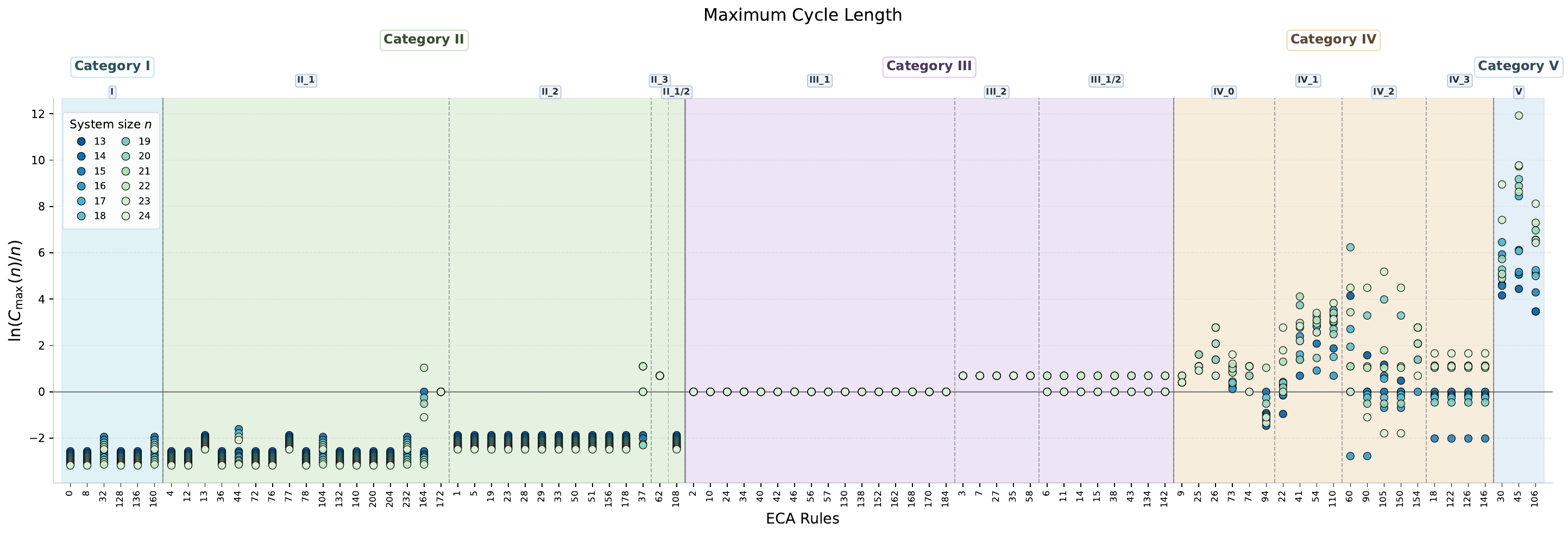}
    \caption{Representation of $\Cmax(n)/n$ for the various ECA in logarithm scale }
    \label{fig_Cmax}
\end{figure}

Since $\Cmax(n)$ plays a role in the definition of Category IV and $V$,
we provide in Figure \ref{fig_Cmax} a plot of $\Cmax(n)/n$ in the logarithmic scale.
For convenience, we also summarize the classification based on spectral data in Table \ref{table_summary},

\begin{table}[htbp]
\centering
\resizebox{\columnwidth}{!}{%
\begin{tabular}{|c|c|l|}
\hline
\textbf{Category}& \textbf{Family} & \textbf{Rule} \\ \hline
I & & 0, 8, 32, 128, 136, 160 \\ \hline
II & II${}_1$ & 4, 12, 13, 36, 44, 72, 76, 77, 78, 104, 132, 140, 164, 172, 200, 204, 232 \\ \hline
& II${}_{1/2}$ & 108 \\ \hline
& II${}_2$ & 1, 5, 19, 23, 28, 29, 33, 37, 50, 51, 156, 178 \\ \hline
& II${}_3$ & 62 \\ \hline
III & III${}_1$ & 2, 10, 24, 34, 40, 42, 46, 56, 57, 130, 138, 152, 162, 168, 170, 184 \\ \hline
& III${}_{1/2}$ & 6, 11, 14, 15, 38, 43, 134, 142 \\ \hline
& III${}_2$ & 3, 7, 27, 35, 58 \\ \hline
IV & IV${}_0$ & 9, 25, 26, 73, 74, 94 \\ \hline
& IV${}_1$ & 22, 41, 54, 110 \\ \hline
& IV${}_2$ & 60, 90, 105, 150, 154 \\ \hline
& IV${}_3$ & 18, 122, 126, 146 \\ \hline
V & & 30, 45, 106 \\ \hline
\end{tabular}
}
\caption{Summary of the classification based on spectral data}
\label{table_summary}
\end{table}

\section{Additional spectral properties}\label{sec_additional}

In this section, we gather various observations about the spectra of Koopman
operators. Most of these observations or results are already known, but studying the
operator $\KRn$ encompasses them in a unique framework.

\subsection{Role of n}

The nature of the value $n\in \N$ with $n\geq 4$ plays an important role for the 
spectral properties of $\KRn$.
As illustrated in the previous section, the parity of $n$ is an important factor. This is not surprising since the 
exchange  of the two elements of the alphabet $\{0,1\}$ is an involution.
In the above classification the role of the parity of $n$ has been emphasized for families I${}_{1/2}$ and III${}_{1/2}$.

\begin{figure}[htbp]
\centering
\begin{minipage}[b]{0.5\textwidth}
\centering
\includegraphics[width=0.9\textwidth]{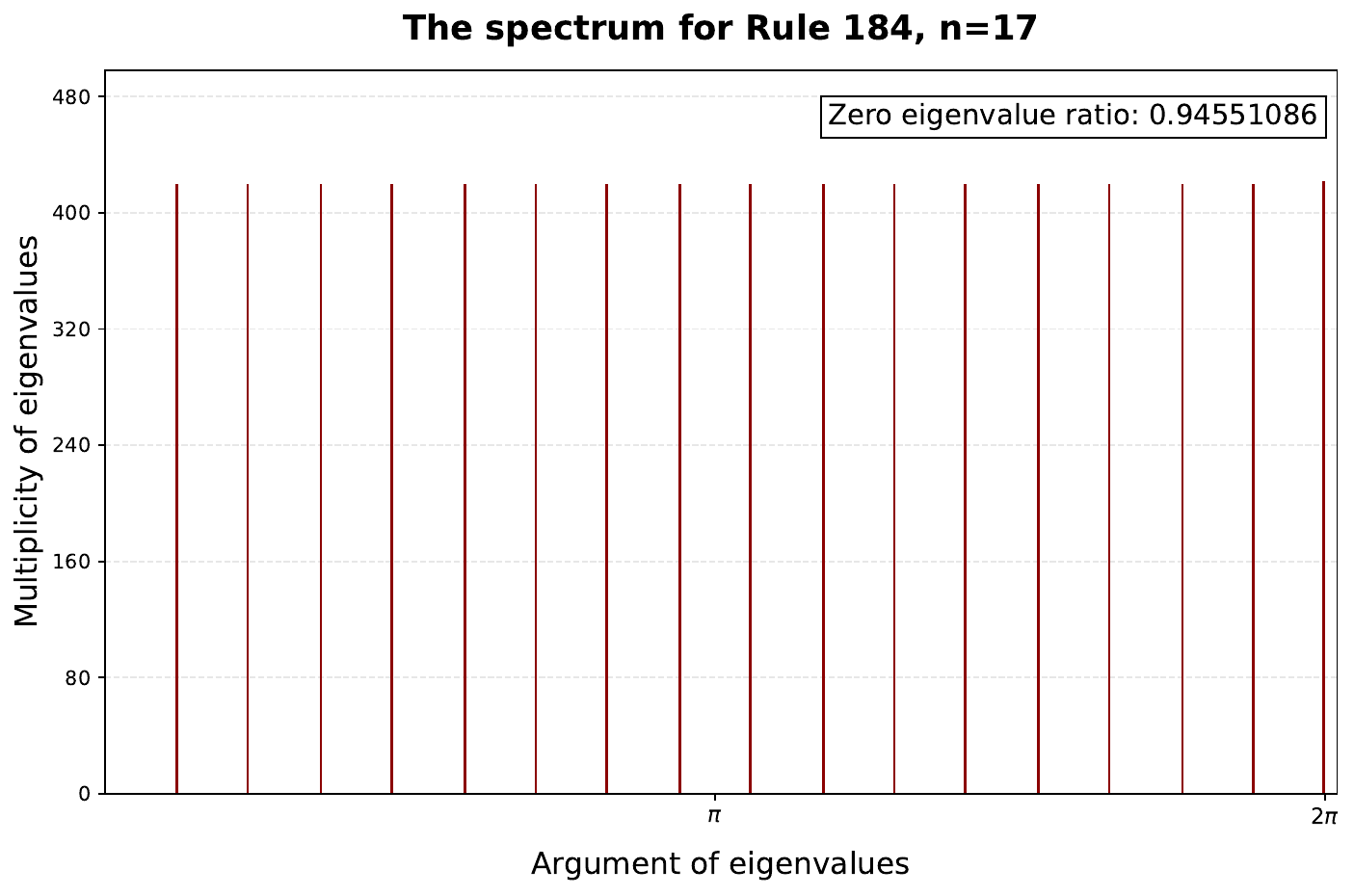}
\end{minipage}%
\begin{minipage}[b]{0.5\textwidth}
\centering
\includegraphics[width=0.9\textwidth]{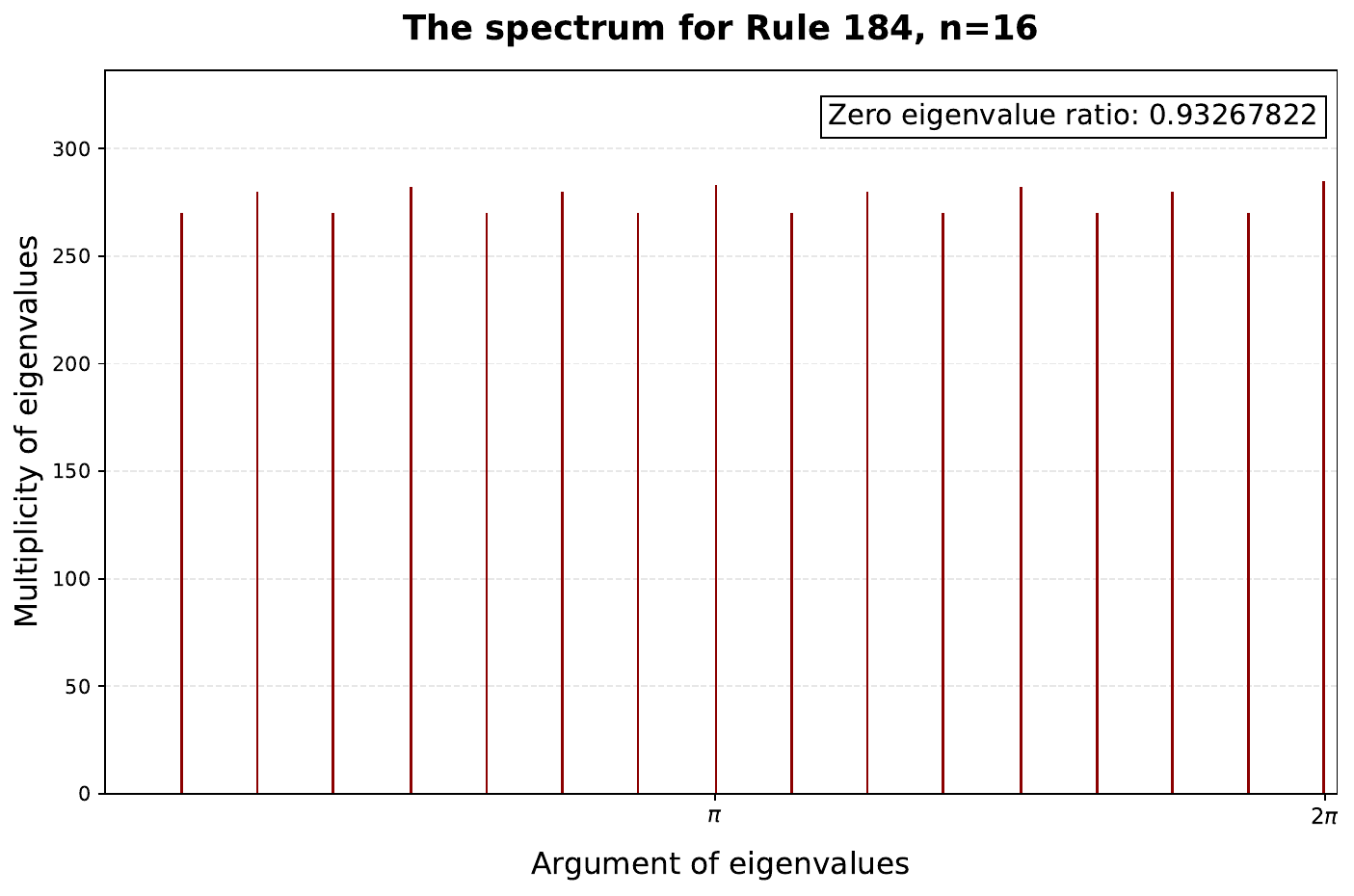}
\end{minipage}%
\caption{Spectra for Rule $184$ for $n=17$ and $n=16$}
\label{fig_prime}
\end{figure}

Another property of $n$ plays an important role, especially for Category III. 
Indeed, for this category the cycle length spectrum $\GRn$ contains the values 
$n$ or $2n$ with large multiplicity, and divisors of these two numbers with smaller multiplicity.  
In the former case, when $n$ is a prime number, then $\GRn$ contains
only the numbers $0$, $1$, and $n$. As a consequence, the corresponding spectrum of $\KRn$ contains only (excluding $0$)
the $n$-th complex root of $1$, all except $1$ of the same multiplicity.
Conversely when $n$ has several divisors, the set $\GRn$ is more complex, and so
is the corresponding spectrum of $\KRn$, see Figure \ref{fig_prime}.
Up to the best of our knowledge, this dependence of $\sigma(\KRn)$ on the nature (prime or not) of $n$ has
not been noticed before.

\subsection{Unitarity}\label{sec_unitary}

As mentioned in Proposition \ref{prop_spectrum} any transient state or any Garden 
of Eden produces an eigenvalue $0$ in the spectrum of $\KRn$. However, there exist
situations for  which none of them exists, which means that $0$ does not belong
to the spectrum of $\KRn$. In these cases the entire spectrum of Koopman
operator is included in the unit circle of the complex plane, and this operator 
is then unitary. Let us mention when this happens.

For Rules $204$ (identity),  $51$ (inversion), $170$ (left shift), and $15$ (right shift with inversion), this is trivially the case. 
More interestingly, for Rules $105$ and $150$, this happens whenever $n$ is not a multiple of $3$, see \cite{{DelRey}}
and earlier references mentioned there for Rule $150$, and
\cite{CLMR,Kr} for relations between Rules $105$ and $150$. 
For Rules $45$ and $154$, Koopman operator is unitary for every odd $n$, as shown in \cite[Table~1]{NobeYura}.
More intriguingly, for $n$ even, Koopman operator $K_{154_n}$ is not far from being unitary, 
as explained below.

\begin{figure}[htbp]
\centering
\begin{minipage}[b]{0.5\textwidth}
\centering
\includegraphics[width=0.9\textwidth]{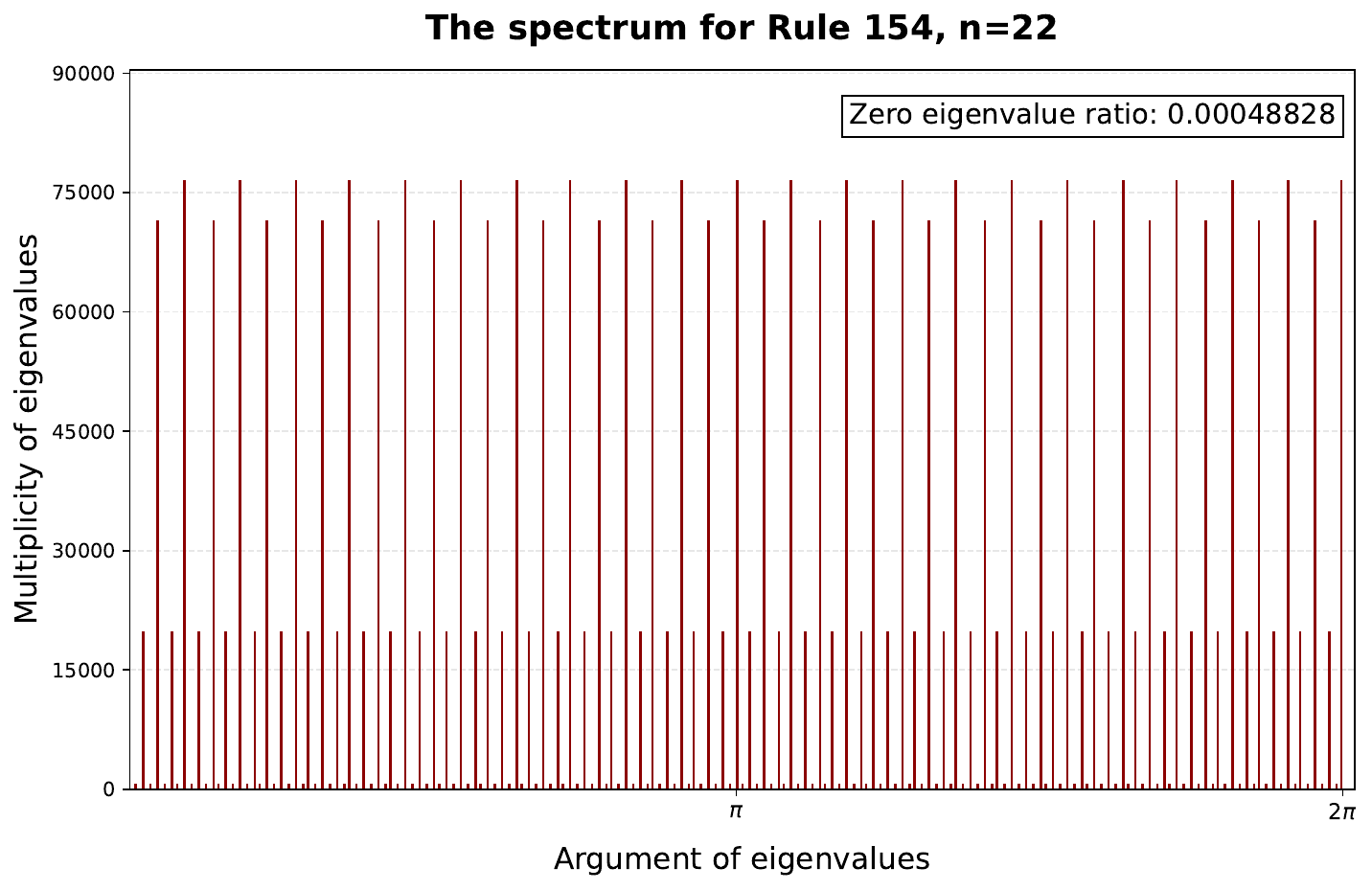}
\end{minipage}%
\begin{minipage}[b]{0.5\textwidth}
\centering
\includegraphics[width=0.9\textwidth]{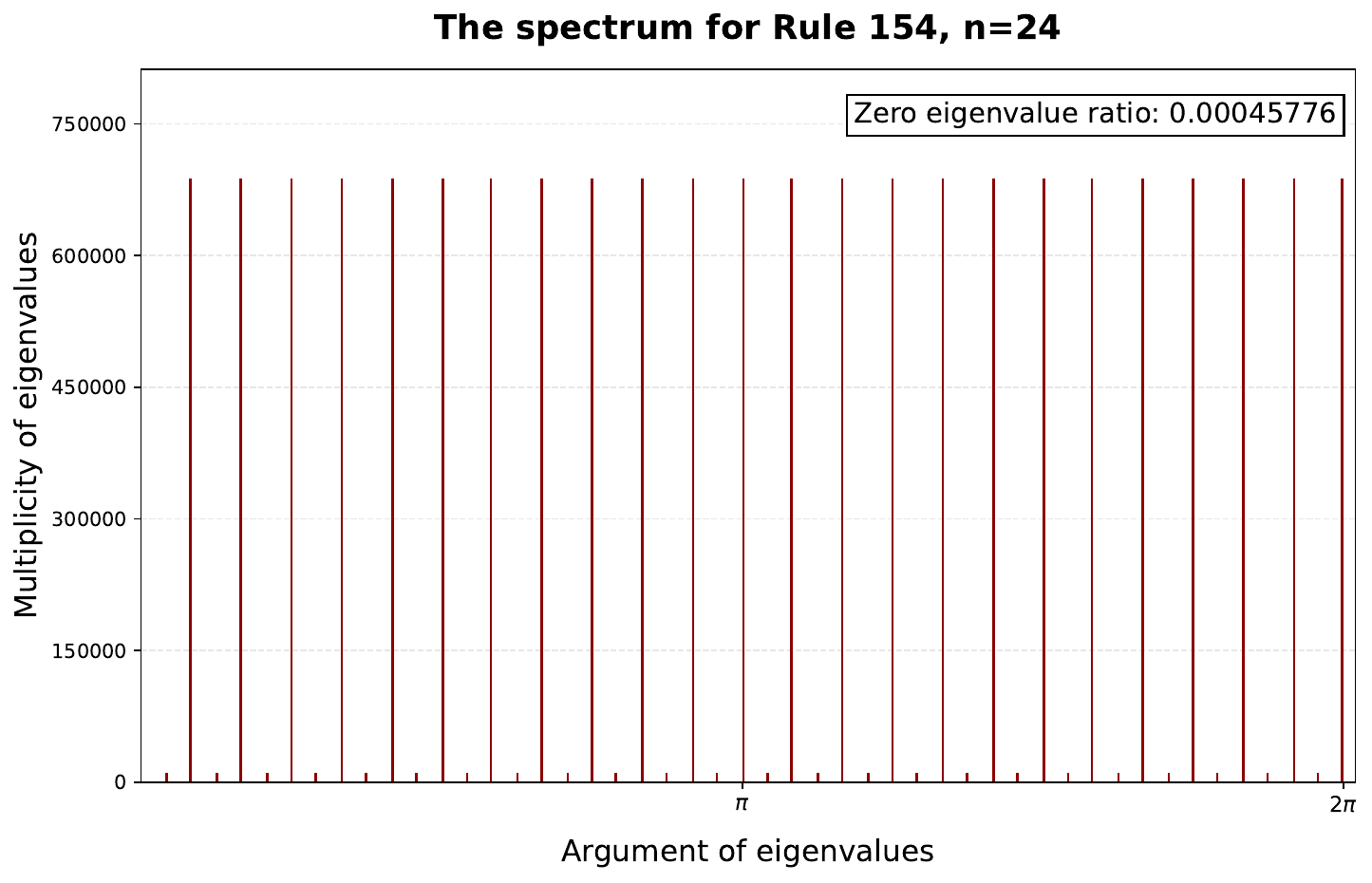}
\end{minipage}%
\caption{Spectra for Rule $154$ for $n=22$ and $n=24$}
\label{fig_154}
\end{figure}

So far, the role of $0$ in the spectrum was essentially overlooked. As already mentioned, we define
the \emph{zero eigenvalue ratio}
\begin{equation}\label{eq_ratio}
\frac{\# \{0\in \sigma(\KRn)\}}{2^n}
\end{equation}
which provides a rough measure of the loss of information of the underlying dynamical system. 
When $\KRn$ is unitary, this ratio is equal to $0$, while very often this ratio is between $0.9$ and $1$. 
However, for Rule $154$ this ratio is always below $0.01$ for $n\in \{13, \dots, 27\}$ and seems to decrease with larger $n$,
see the zero eigenvalue ratio indicated in Figure \ref{fig_154}.
For comparison, $K_{45_n}$ is also unitary for $n$ odd, 
but for $n$ even and belonging to $\{13, \dots, 24\}$ the zero eigenvalue ratio is always bigger than $0.5$.

\subsection{Families in Category IV}

Category IV contains various ECA which do not fit in other categories and for which 
we could not find any specific pattern in the spectra of $\KRn$ for different $n$.
Inside this category, some families of ECA can now be defined.

The most interconnected family is made of Rules $18$, $122$, $126$ and $146$. Very surprisingly, the spectra of
the corresponding Koopman operators for these four ECA are very similar, for all $n\in \{13, \dots, 27\}$. 
In fact, despite having many eigenvalues (growing with $n$), their spectra differ at most by 
a change of multiplicity of the eigenvalues $\pm 1$.
At the level of the cycle length spectrum $\GRn$, the only difference is the appearance
of one additional fixed state $C_1$ or of one additional periodic orbit $C_2$.
From a different perspective these similarities had already been observed visually in \cite{matho}
and in \cite{Jiang}.

Another family is made of Rules $60$, $90$, $105$, $150$, and $154$. Firstly, they share the common feature
of being balanced (meaning their code in binary numbers contains four number $1$ and four numbers $0$).
Most importantly, the ratio providing the information of the loss of information introduced in \eqref{eq_ratio}
has regular patterns that are not shared by other ECA. For Rules $105$, $150$ and $154$ these patterns
have already been mention in Section \ref{sec_unitary}. On the other hand, Rules $60$ and $90$
do not lead to a unitary operator $\KRn$ (meaning the absence of eigenvalue $0$), but their zero eigenvalue ratio are $0.5$ 
for odd $n\in \{13, \dots, 24\}$, $0.75$ for $n\in \{14,18,22\}$, and above 
$0.93$ for $n\in \{16, 20, 24\}$. We do not have an explanation for this surprising pattern.

For the remaining ECA in Category IV, any additional subdivision is less clear. The only additional characteristic
we have found which could provide two families is the size of the largest cycle of the cycle length spectrum $\GRn$.
However, precisely because the values inside the cycle length spectrum $\GRn$ are quite difficult to predict
for different $n$ for the ECA in this category, we do not propose any further subdivisions.

\section{Epilogue}

The various results presented above have all been obtained by looking at finite periodic realizations
of ECA obtained for $n$ in a sufficiently large set. 
The next step is to investigate the infinitely extended realizations, when ECA are considered on $\Z$.
For that purpose, tracking the spectral distributions as the system sizes increase will be essential. 
In particular, it will be necessary to investigate which spectral features remain stable and 
which are caused mainly by the finite size of the systems. For properties
which highly depend on $n$, as for example some mentioned in Section \ref{sec_unitary}, 
it is not clear what is going to happen. We intend to investigate these questions in the future.

\end{document}